\documentclass[aap,MSNbibl,seceqn,citesort,dvips]{arximspdf}

\doi{10.1214/11-AAP762}
\volume{21}
\issue{6}
\pubyear{2011}
\firstpage{2379}
\lastpage{2423}

\makeatletter

\newtheorem{theorem}{Theorem}[section]
\newtheorem{corollary}[theorem]{Corollary}
\newtheorem{lemma}[theorem]{Lemma}

\newproclaim{definition}[theorem]{Definition}
\newproclaim{example}[theorem]{Example}
\newproclaim{remark}[theorem]{Remark}
\newproclaim{assumption}{Assumption}[section]

\newcommand{\EE}{\mathbb{E}}

\newcommand{\cF}{\mathcal{F}}

\makeatother

\begin{document}
\begin{frontmatter}

\title{Malliavin calculus for backward stochastic differential
equations and application to numerical solutions}
\runtitle{Malliavin calculus, numerical solution of bsde}

\begin{aug}
\author[A]{\fnms{Yaozhong} \snm{Hu}\corref{}\thanksref{t1}\ead[label=e1]{hu@math.ku.edu}},
\author[A]{\fnms{David} \snm{Nualart}\thanksref{t2}\ead[label=e2]{nualart@math.ku.edu}} and
\author[A]{\fnms{Xiaoming} \snm{Song}\ead[label=e3]{xsong@math.ku.edu}}
\runauthor{Y. Hu, D. Nualart and X. Song}
\affiliation{University of Kansas}
\address[A]{Department of Mathematics\\
University of Kansas \\
Lawrence, Kansas 66045\\
USA\\
\printead{e1}\\
\hphantom{E-mail: }\printead*{e2}\\
\hphantom{E-mail: }\printead*{e3}} %adresu isvedimo komanda gale!
\end{aug}

\thankstext{t1}{Supported by the NSF Grant DMS-05-04783.}
\thankstext{t2}{Supported by the NSF grant DMS-09-04538.}

% HISTORY:
\received{\smonth{3} \syear{2010}}
\revised{\smonth{10} \syear{2010}}

% ABSTRACT
%
\begin{abstract}
In this paper we study backward stochastic differential equations
with general terminal value and general random generator. In
particular, we do not require the terminal value be given by a
forward diffusion equation. The randomness of the generator does
not need to be from a forward equation, either. Motivated from
applications to numerical simulations, first we obtain the
$L^p$-H\"{o}lder continuity of the solution. Then we construct
several numerical approximation schemes for backward stochastic
differential equations and obtain the rate of convergence of the
schemes based on the obtained $L^p$-H\"{o}lder continuity results.
The main tool is the Malliavin calculus.
\end{abstract}

% KEYWORDS
%

\pdfkeywords{60H07, 60H10, 60H35, 65C30, 91G60, Backward stochastic differential equations,
Malliavin calculus, explicit scheme, implicit scheme,
Clark--Ocone--Haussman formula, rate of convergence,
Holder continuity of the solutions}

\end{frontmatter}

%s1 ###
\section{Introduction}

The backward stochastic
differential equation (BSDE, for short) we shall consider in this paper
takes the following form:
%
%e1.1 ###
%
\begin{equation} \label{bsde}
Y_t=\xi+\int_t^Tf(r, Y_r,Z_r)\,dr-\int_t^TZ_r\,dW_r,\qquad 0\le t\le T ,
\end{equation}
where $W=\{W_t\}_{0\le t\le T}$ is a standard Brownian motion,
$\xi$ is the given terminal value and $f$ is the given (random)
generator. To solve this equation is to find a pair of adapted
processes $Y=\{Y_t\}_{0\le t\le T}$ and $Z=\{Z_t\}_{0\le t\le T}$
satisfying the above equation (\ref{bsde}).

Linear backward stochastic differential equations were first
studied by Bismut \cite{bismut} in an attempt to solve some
optimal stochastic control problem through the method of maximum
principle. The general nonlinear backward stochastic differential
equations were first studied by Pardoux and Peng \cite{PP90}.
Since then there have been extensive studies of this equation. We
refer to the review paper
by El Karoui,
Peng and Quenez \cite{KPQ97}, to the books of El Karoui and
Mazliak \cite{KM} and of Ma and Yong \cite{MY} and the references therein
for more comprehensive presentation of the theory.

A current important topic in the applications of BSDEs is the numerical
approximation schemes. In most work on numerical simulations,\vadjust{\goodbreak}
a certain forward stochastic differential equation of the following form:
%
%e1.2 ###
%
\begin{equation} \label{e.1.2}
X_t=X_0+\int_0^tb(r, X_r, Y_r)\,dr+\int_0^t{\sigma}(r, X_r)\,dW_r
\end{equation}
is needed. Usually it is assumed that
the generator $f$ in (\ref{bsde}) depends on~$X_r$ at the
time $r$: $f(r, Y_r, Z_r)=f(r, X_r, Y_r, Z_r)$, where $f(r, x, y,
z)$ is a deterministic function of $(r, x, y, z)$, and $f$ is
global Lipschitz in $(x,y,z)$. If in addition the terminal value
$\xi$ is of the form $ \xi=h(X_T)$, where $h$ is a~deterministic
function, a~so-called four-step numerical scheme has been
developed by Ma, Protter and Yong in \cite{MPY}. A basic
ingredient in this paper is that the solution $\{Y_t\}_{0\le t\le
T}$ to the BSDE is of the form $Y_t=u(t, X_t)$, where $u(t,x)$ is
determined by a quasi-linear partial differential equation of
parabolic type. Recently, Bouchard and Touzi \cite{BT} propose a
Monte-Carlo approach which may be more suitable for
high-dimensional problems. Again in this forward--backward setting,
if the generator $f$ has a quadratic
growth in $Z$, a numerical approximation is developed by Imkeller
and Dos Reis \cite{IDR} in which a truncation procedure is
applied.

In the case where the terminal value $\xi$ is a functional of the
path of the forward diffusion $X$, namely, $ \xi=g(X_\cdot)$,
different approaches to construct numerical methods have been
proposed. We refer to Bally \cite{Ba} for a scheme with a random
time partition. In the work by Zhang \cite{Zh}, the
$L^2$-regularity of $Z$ is obtained, which allows one to use
deterministic time partitions as well as to obtain the rate
estimate (see Bender and Denk \cite{BD}, Gobet, Lemor and Warin
\cite{GLW} and Zhang \cite{Zh} for different algorithms). We
should also mention the works by Briand, Delyon and M\'{e}min
\cite{BDM} and Ma et al. \cite{MPSS},
where the Brownian motion is replaced by a scaled random walk.

The purpose of the present paper is to construct numerical
schemes for the general
BSDE (\ref{bsde}), without assuming any particular form for the
terminal value $\xi$ and generator $f$. This means that $\xi$
can be an arbitrary random variable, and $f(r, y, z)$ can be an
arbitrary $\cF_r$-measurable random variable (see Assumption \ref{a.3.2}
in Section~\ref{sec2} for precise conditions on $\xi$ and $f$). The
natural tool that we shall use is the Malliavin calculus. We
emphasize that the main difficulty in constructing a numerical
scheme for BSDEs is usually the approximation of the process
$Z$. It is necessary to obtain some regularity properties for
the trajectories of this process~$Z$. The Malliavin calculus turns
out to be a suitable tool to handle these problems because the
random variable $Z_t$ can be expressed in terms of the trace of
the Malliavin derivative of $Y_t$, namely, $Z_t=D_tY_t$. This
relationship was proved in the paper by El Karoui, Peng and
Quenez \cite{KPQ97} and was used by these authors to obtain estimates
for the moments of $Z_t$. We shall further exploit this identity
to obtain the $L^p$-H\"{o}lder continuity of the process $Z$,
which is the critical ingredient for the rate estimate of our
numerical schemes.

Our first numerical scheme was inspired by the paper of
Zhang \cite{Zh}, where the author considers a class of BSDEs whose
terminal value $\xi$ takes the form~$g(X_{\cdot})$, where $X$ is
a forward diffusion of the form (\ref{e.1.2}), and $g$ satisfies
a~Lipschitz condition with respect to the $L^\infty$ or $L^1$
norms (similar assumptions for $f$). The discretization scheme is
based on the regularity of the process $Z$ in the mean square
sense; that is,
for any partition $\pi=\{0=t_0<t_1<\cdots<t_n=T\}$, one obtains
%
%e1.3 ###
%
\begin{equation}\label{e.1.3}
\sum_{i=0}^{n-1} \mathbb{E} \int_{t_i}^{t_{i+1}} [
|Z_t-Z_{t_i}|^2+ |Z_t-Z_{t_{i+1}}|^2] \,dt\le K |\pi| ,
\end{equation}
where $|\pi|=\max_{0\le i\le n-1} (t_{i+1}-t_i)$, and $K$ is a
constant independent of the partition $\pi$.

We consider the case of a general terminal value $\xi$ which is
twice differentiable in the sense of Malliavin calculus, and the
first and second derivatives satisfy some integrability conditions; we
also made similar assumptions for the generator $f$
(see Assumption \ref{a.3.2} in Section \ref{sec2} for details). In this
sense our framework
extends that of \cite{MZh} and is also natural. In this framework, we are
able to obtain an estimate of the form
%
%e1.4 ###
%
\begin{equation} \label{e.1.4}
\mathbb{E} |Z_t-Z_s|^p\le K|t-s|^{p/2} ,
\end{equation}
where $K$ is a constant independent of $s$ and $t$. Clearly,
(\ref{e.1.4}) with $p=2$ implies (\ref{e.1.3}). Moreover,
(\ref{e.1.4}) implies the existence of a $\gamma$-H\"{o}lder
continuous version of the process $Z$ for any $\gamma<\frac
{1}{2}-\frac{1}{p}$. Notice that, up to now the path regularity of
$Z$ has been studied only when the terminal value
and the generator are
functional of a forward diffusion.

After establishing the regularity of $Z$, we consider different
types of numerical schemes. First we analyze a scheme similar to the
one proposed
in~\cite{Zh} [see (\ref{e.4.2})]. In this case we obtain a
rate of convergence of the following type:
\[
\mathbb{E} \sup_{0\le t\le T} |Y_{t}-Y_{t}^\pi|^2+\int_{0}^T
\mathbb{E} |Z_t-Z_t^\pi|^2 \,dt\le
K(|\pi|+\EE|\xi-\xi^\pi|^2) .
\]
Notice that this result is stronger than that in \cite{Zh} which can
be stated as
(when $\xi^\pi=\xi$)
\[
\sup_{0\le t\le T} \mathbb{E} |Y_{t}-Y_{t}^\pi|^2+\int_{0}^T
\mathbb{E} |Z_t-Z_t^\pi|^2 \,dt\le
K |\pi| .
\]

We also propose and study an ``implicit'' numerical scheme [see %
(\ref{e.5.1}) in Section~\ref{sec4} for the details]. For this
scheme we obtain a much better result on the rate of convergence,
\[
\mathbb{E} \sup_{0\le t\le T}
|Y_{t}-Y_{t}^\pi|^p+\mathbb{E}\biggl(\int_{0}^T
|Z_t-Z_t^\pi|^2 \,dt\biggr)^{p/2}\le K(|\pi|^{p/2}+\EE
|\xi-\xi^\pi|^p) ,
\]
where $p>1$ depends on the assumptions imposed on the terminal value
and the coefficients.

In both schemes, the integral of the process $Z$ is used in each
iteration, and for this reason they are not completely discrete
schemes. In order to implement the scheme on computers, one must
replace an integral of the form $\int_{t_{i }}^{t_{i+1}} Z_s^\pi
\,ds $ by discrete sums, and then the convergence of the obtained
scheme is hardly guaranteed. To avoid this discretization we
propose a truly discrete numerical scheme using our representation
of $Z_t$ as the trace of the Malliavin derivative of $Y_t$ (see
Section~\ref{sec5} for details). For this new scheme, we obtain a rate of
convergence result of the form
\[
\mathbb{E} \max_{0\le i\le n} \{ |Y_{t_i}-Y_{t_i}^\pi|^p+
|Z_{t_i}-Z_{t_i}^\pi|^p \} \le K|\pi|^{
p/2-{\varepsilon}}
\]
for any ${\varepsilon}>0$. In fact, we have a slightly better rate
of convergence (see Theorem~\ref{t.6.1}),
\[
\mathbb{E} \max_{0\le i\le n} \{ |Y_{t_i}-Y_{t_i}^\pi|^p +
|Z_{t_i}-Z_{t_i}^\pi|^p\}\le K |\pi|^{{p/2}-
{p}/({2\log({1}/{|\pi|})%
})} \biggl( \log\frac{1}{|\pi|}\biggr)^{p/2}.
\]
However, this type of result on the rate of convergence applies
only to some classes of BSDEs, and thus this scheme remains to be
further investigated.

In the computer realization of our schemes or any other schemes,
an extremely important procedure is to compute the conditional
expectation of form $\EE(Y|\cF_{t_i})$. In this paper we shall
not discuss this issue but only mention the papers \cite{BD,BT} and~\cite{GLW}.

The paper is organized as follows. In Section \ref{sec2} we obtain a
representation of the martingale integrand $Z$ in terms of the
trace of the Malliavin derivative of $Y$, and then we get the
$L^p$-H\"{o}lder continuity of $Z$ by using this representation.
The conditions that we assume on the terminal value~$\xi$ and the
generator $f$
are also specified in this section. Some examples of application are presented
to explain the validity of the conditions. Section~%
\ref{sec3} is devoted to the analysis of the approximation scheme
similar to the one introduced in
\cite{Zh}. Under some differentiability and
integrability conditions in the sense of Malliavin calculus on
$\xi$ and the nonlinear coefficient $f$, we establish a better
rate of convergence for this scheme. In Section~\ref{sec4}, we introduce
an ``implicit'' scheme and obtain the rate of convergence in the
$L^p$ norm. A~completely discrete scheme is proposed and analyzed
in Section~\ref{sec5}.

Throughout the paper for simplicity we consider only scalar BSDEs. The
results obtained in this paper can be easily extended to
multi-dimensional BSDEs.

%s2 ###
\section{The Malliavin calculus for BSDEs}\label{sec2}

%s2.1 ###
\subsection{Notations and preliminaries}

Let $W=\{ W_{t}\}_{0\le t\le T} $ be a
one-dimen\-sional standard Brownian motion defined on some complete
filtered probability space $(\Omega
, \mathcal{F}, P, \{\mathcal{F}_{t}\}_{0\leq t\leq T})$.
We assume that $\{\mathcal{F}_{t}\}_{0\le t\le T} $ is the
filtration generated by the Brownian motion and the $P$-null sets, and
$\mathcal{F}=%
\mathcal{F}_{T}$. We denote by $\mathcal{P}$ the progressive $\sigma
$-field on the
product space $[0,T]\times\Omega$.

For any $p\geq1$ we consider
the following classes of processes:
\begin{itemize}
\item$M^{2,p}$, for any $p\geq2$, denotes the class of square
integrable random variables $F$
with a stochastic integral representation of the form%
\[
F=\mathbb{E}F+\int_{0}^{T}u_{t}\,dW_{t},
\]
where $u$ is a progressively measurable process satisfying $\sup
_{0\leq t\leq T}\mathbb{E}%
|u_{t}|^{p}<\infty$.
\item$H_{\mathcal{F}}^{p}([0,T])$ denotes the Banach space
of all
progressively measurable processes $\varphi\dvtx([0,T]\times\Omega%
, \mathcal{P})\rightarrow(\mathbb{R}, \mathcal{B})$ with norm%
\[
\Vert\varphi\Vert_{H^{p}}=\biggl( \mathbb{E}\biggl( \int
_{0}^{T}|\varphi
_{t}|^{2}\,dt\biggr) ^{{p/2}}\biggr) ^{{1/p}}<\infty.
\]

\item$S_{\mathcal{F}}^{p}([0,T])$ denotes the Banach space of all
the RCLL (right continuous with left limits) adapted processes
$\varphi\dvtx([0,T]\times
\Omega, \mathcal{P})\rightarrow
(\mathbb{R}, %
\mathcal{B})$ with norm
\[
\Vert\varphi\Vert_{S^{p}}=\Bigl(\mathbb{E} \sup_{0\leq t\leq
T}|\varphi_{t}|^{p} \Bigr)^{{1/p}}<\infty.
\]
\end{itemize}

Next, we present some preliminaries on Malliavin calculus, and we
refer the reader to the book by Nualart \cite{N06} for
more details.

Let $\mathbf{H}=L^2([0,T])$ be the separable Hilbert space of all
square integrable real-valued functions on the interval $[0,T]$
with scalar product denoted by $\langle
\cdot,\cdot\rangle_\mathbf{H}$. The norm of an element $h\in
\mathbf{H}$ will be denoted by $\Vert h\Vert_\mathbf{H}$. For any
$h\in\mathbf{H}$ we put $W(h)=\int_0^Th(t)\,dW_t$.

We denote by $C_p^\infty(\mathbb{R}^n)$ the set of all infinitely
continuously differentiable functions $g\dvtx\mathbb{R}^n\rightarrow
\mathbb{R}$ such that $g$ and all of its partial derivatives have
polynomial growth. We make use of the notation $\partial_i g=\frac
{\partial
g}{\partial x_i}$ whenever $g\in C^1(\mathbb{R}^n)$.

Let $\mathcal{S}$ denote the class of smooth random variables such
that a random variable $F\in\mathcal{S}$ has the form
%
%e2.1 ###
%
\begin{equation}\label{smooth}
F=g(W(h_1),\ldots,W(h_n)),
\end{equation}
where $g$ belongs to $C_p^\infty(\mathbb{R}^n)$, $h_1,\ldots,h_n$
are in $\mathbf{H}$ and $n\geq1$.

The Malliavin derivative of a smooth random variable $F$ of the
form (\ref{smooth}) is the $\mathbf{H}$-valued random variable
given by
\[
DF=\sum_{i=1}^n\partial_i g(W(h_1),\ldots,W(h_n))h_i.\vadjust{\goodbreak}
\]
For any $p\geq1$ we will denote the domain of $D$ in $L^p(\Omega)$
by $\mathbb{D}^{1,p}$, meaning that~$\mathbb{D}^{1,p}$ is the
closure of the class of smooth random variables $\mathcal{S}$ with
respect to the norm
\[
\Vert F\Vert_{1,p}=(\EE|F|^p+\EE\Vert
DF\Vert_\mathbf{H}^p)^{1/p}.
\]
We can define the iteration of the operator $D$ in such a way that
for a smooth random variable $F$, the iterated derivative $D^kF$
is a random variable with values in $\mathbf{H}^{\otimes k}$. Then
for every $p\geq1$ and any natural number $k\geq1$ we introduce
the seminorm on $\mathcal{S}$ defined by
\[
\Vert F\Vert_{k,p}=\Biggl(\EE|F|^p+\sum_{j=1}^k\EE\Vert
D^jF\Vert_{\mathbf{H}^{\otimes j}}^p\Biggr)^{1/p}.
\]
We will denote by $\mathbb{D}^{k,p}$ the completion of the family
of smooth random variables $\mathcal{S}$ with respect to the norm
\mbox{$\Vert\cdot\Vert_{k,p}$}.

Let $\mu$ be the Lebesgue measure on $[0,T]$. For any $k\geq1$
and $F\in\mathbb{D}^{k,p}$, the derivative
\[
D^kF=\{D^k_{t_1,\ldots,t_k}F, t_i\in[0,T], i=1,\ldots,k\}
\]
is a measurable function on the product space
$[0,T]^k\times\Omega$, which is defined a.e. with respect to the
measure $\mu^k\times P$.

We use $\mathbb{L}_{a}^{1,p}$ to denote the set of
real-valued progressively measurable processes $u=\{
u_t\}_{0\leq t\leq T} $ such that:

\begin{longlist}
\item For almost all $t\in\lbrack0,T], u_t\in
\mathbb{D}^{1,p}$.\vspace*{2pt}

\item$\mathbb{E}(
(\int_{0}^{T}|u_t|^{2}\,dt)^{{p/2}}+(\int
_{0}^{T}\int_{0}^{T}|D_{\theta
}u_t|^{2}\,d\theta \,dt)^{{p/2}}) <\infty$.
\end{longlist}
Notice that we can choose a progressively measurable version of the $%
\mathbf{H}$-valued process $\{Du_{t}\}_{0\le t\le T}$.

%s2.2 ###
\subsection{Estimates on the solutions of BSDEs}
The generator $f$ in the BSDE (\ref{bsde}) is a measurable function
$f\dvtx([0,T]\times
\Omega\times\mathbb{R}\times\mathbb{R}, \mathcal{P}\otimes
\mathcal{B}\otimes\mathcal{B})\rightarrow
(\mathbb{R}, \mathcal{B})$, and the terminal value $\xi$ is an
$\mathcal{F}_T$-measurable
random variable.
\begin{definition}
A solution to the BSDE (\ref{bsde}) is a pair of progressively
measurable processes $(Y,Z)$ such that $\int_0^T|Z_t|^2\,dt\,{<}\,\infty$,
\mbox{$\int_0^T|f(t,Y_t,Z_t)|\,dt\,{<}\,\infty$}, a.s. and
\[
Y_t=\xi+\int_t^Tf(r, Y_r,Z_r)\,dr-\int_t^TZ_r\,dW_r,\qquad 0\le t\le T.
\]
\end{definition}

The next lemma provides a useful estimate on the solution to the
BSDE~(\ref{bsde}).
\begin{lemma}
\label{T.2.1} Fix $q\geq2$. Suppose that $\xi\in L^{q}(\Omega)$,
$f(t,0,0)\in H_{\mathcal{F}}^{q}([0,T])$ and $f$ is uniformly
Lipschitz in $(y,z)$; namely, there exists a positive number~$L$
such that $\mu\times P $ a.e.
\[
|f(t,y_{1},z_{1})-f(t,y_{2},z_{2})|\leq
L(|y_{1}-y_{2}|+|z_{1}-z_{2}|)
\]
for all $y_1,y_2\in\mathbb{R}$ and $z_1,z_2\in\mathbb{R}$.
Then
there exists a unique solution pair $(Y,Z)\in
S_{\mathcal{F}}^{q}([0,T])\times H_{\mathcal{F}}^{q}([0,T])$ to
(\ref{bsde}). Moreover, we have the following estimate
for the solution:
%
%e2.2 ###
%
\begin{eqnarray} \label{e.2.1}
&&
\mathbb{E}\sup_{0\leq t\leq T}|Y_{t}|^{q}+\mathbb{E}\biggl(
\int_{0}^{T}|Z_{t}|^{2}\,dt\biggr) ^{{q/2}}\nonumber\\[-8pt]\\[-8pt]
&&\qquad\leq K\biggl(
\mathbb{E}|\xi
|^{q}+\mathbb{E}\biggl( \int_{0}^{T}|f(t,0,0)|^{2}\,dt\biggr) ^{
{q/2}}\biggr) ,\nonumber
\end{eqnarray}
where $K$ is a constant depending only on $L$, $q$ and $T$.
\end{lemma}
\begin{pf}
The proof of the existence and uniqueness of the solution $(Y,Z)$
can be found in \cite{KPQ97}, Theorem 5.1, with the local
martingale $M\equiv0$, since the filtration here is the filtration
generated by the Brownian motion $W$. Estimate (\ref{e.2.1})
can be easily obtained from Proposition 5.1 in \cite{KPQ97} with
$(f^1,\xi^1)=(f,\xi)$ and $(f^2,\xi^2)=(0,0)$.
\end{pf}

As we will see later, for a given BSDE the process $Z$ will be
expressed in terms of the Malliavin derivative of the solution $Y$, which
will satisfy a~linear BSDE with random coefficients. To study the
properties of $Z$ we need to analyze a class of linear BSDEs.

Let $\{{\alpha}_{t}\}_{0\le t\le T}$ and $\{{\beta}_{t}\}_{0\le
t\le T}$ be two progressively measurable processes. We will make
use of the following integrability conditions:
\begin{assumption}\label{a.2.1}
\begin{longlist}[(H1)]
\item[(H1)]
For any $\lambda>0$,
\[
C_{\lambda}:=\mathbb{E} \exp\biggl( \lambda\int_{0}^{T}(
| \alpha_{t}| +\beta_{t}^{2}) \,dt\biggr) <\infty.
\]

\item[(H2)] For any $p\geq1$,
\[
K_{p}:=\sup_{0\leq t\leq T}\mathbb{E}( |\alpha
_{t}|^{p}+|\beta_{t}|^{p}) <\infty.
\]
\end{longlist}
\end{assumption}

Under condition (H1), we denote by $\{\rho_t\}_{0\le t\le T}$ the
solution of the linear stochastic differential equation
%
%e2.3 ###
%
\begin{equation}\label{joint}
\cases{d\rho_{t}=\alpha_{t}\rho_{t}\,dt+\beta_{t}\rho_{t}\,dW_{t}, &\quad
$0\le t\le T$,
\cr
\rho_{0}=1.}
\end{equation}

The following theorem is a critical tool for the proof of the main
theorem in this section, and it has also its own interest.\vadjust{\goodbreak}
\begin{theorem}
\label{T.3.2} Let $q>p\geq2$ and let $\xi\in L^q(\Omega)$ and $%
f\in H_{\mathcal{F}}^{q}([0,T])$. Assume that $\{{\alpha}_{t}\}_{0\le
t\le T}$ and $\{{%
\beta}_{t}\}_{0\le t\le T}$ are two progressively measurable
processes satisfying conditions \textup{(H1)} and \textup{(H2)} in Assumption
\ref{a.2.1}. Suppose that the random variables $\xi\rho_{T}$ and
$\int_{0}^{T}\rho_{t}f_{t}\,dt$ belong to $M^{2,q}$, where $\{\rho_t\}
_{0\le t\le T}$ is the solution to (\ref{joint}).
Then the following linear BSDE,
%
%e2.4 ###
%
\begin{equation}\label{linear}\quad
Y_{t}=\xi+\int_{t}^{T}[\alpha_{r}Y_{r}+\beta
_{r}Z_{r}+f_{r}]\,dr-\int_{t}^{T}Z_{r}\,dW_{r},\qquad 0\le t\le T,
\end{equation}
has a unique solution pair $(Y, Z)$, and there is
a constant $K>0$ such that
%
%e2.5 ###
%
\begin{equation}\label{e.2.4}
\mathbb{E}|Y_{t}-Y_{s}|^{p}\leq K|t-s|^{{p/2}}
\end{equation}
for all $s, t\in[0,T]$.
\end{theorem}

We need the following lemma to prove the above result.
\begin{lemma}
\label{l.3.1} Let $\{{\alpha}_{t}\}_{0\le t\le T}$ and $\{{\beta}%
_{t}\}_{0\le t\le T}$ be two progressively measurable processes
satisfying condition \textup{(H1)} in Assumption \ref{a.2.1}, and $\{\rho
_{t}\}_{0\le t\le T}$ be the solution of (\ref{joint}).
Then, for any $r\in\mathbb{R}$ we have
%
%e2.6 ###
%
\begin{equation}\label{e.3.3}
\mathbb{E} \sup_{0\leq t\leq T}\rho_{t}^{r} <\infty.
\end{equation}
\end{lemma}
\begin{pf}
Let $t\in[0,T]$. The solution to (\ref{joint}) can be
written as
\[
\rho_{t}=\exp\biggl\{ \int_{0}^{t}\biggl( \alpha_{s}-\frac{\beta
_{s}^{2}}{2%
}\biggr)\, ds+\int_{0}^{t}\beta_{s}\,dW_{s}\biggr\} .
\]
For any real number $r$, we have
\begin{eqnarray*}
\mathbb{E} \sup_{0\leq t\leq T}\rho_{t}^{r} &=&\mathbb{E}%
\sup_{0\leq t\leq T}\exp\biggl\{ \int_{0}^{t}r\biggl( \alpha_{s}-%
\frac{\beta_{s}^{2}}{2}\biggr) \,ds+r\int_{0}^{t}\beta
_{s}\,dW_{s}\biggr\} \\
&\leq&\mathbb{E} \biggl( \exp\biggl\{ | r|
\int_{0}^{T}| \alpha_{s}|
\,ds+\frac{1}{2}(|r|+r^{2})\int_{0}^{T}\beta
_{s}^{2}\,ds\biggr\} \\
&&\hspace*{21.2pt}{} \times\sup_{0\leq t\leq T}\exp\biggl\{
r\int_{0}^{t}\beta_{s}\,dW_{s}-\frac{r^{2}}{2}\int_{0}^{t}\beta
_{s}^{2}\,ds\biggr\} \biggr) .
\end{eqnarray*}
Then, fixing any $p>1$ and using H\"{o}lder's inequality, we obtain%
%
%e2.7 ###
%
\begin{equation}\label{martingale}\qquad
\mathbb{E} \sup_{0\leq t\leq T}\rho_{t}^{r} \leq C\biggl(
\EE\sup_{0\leq t\leq T}\exp\biggl\{ rp\int_{0}^{t}\beta
_{s}\,dW_{s}-\frac{%
pr^2}{2} \int_{0}^{t}\beta_{s}^{2}\,ds\biggr\}\biggr) ^{{1}/{p}%
},
\end{equation}
where
\[
C=\biggl( \mathbb{E} \exp\biggl\{ q| r|
\int_{0}^{T}| \alpha
_{s}| \,ds+\frac{q}{2}(|r|+r^{2})\int_{0}^{T}\beta
_{s}^{2}\,ds\biggr\}
\biggr) ^{1/q}
\]
and $\frac{1}{p}+\frac{1}{q}=1$.\vadjust{\goodbreak}

Set $M_{t}=\exp\{ r\int_{0}^{t}\beta_{s}\,dW_{s}-\frac{%
r^{2}}{2}\int_{0}^{t}\beta_{s}^{2}\,ds\}$. Then
$\{M_{t}\}_{0\le t\le T} $ is a martingale
due to (H1). We can rewrite (\ref{martingale}) into
%
%e2.8 ###
%
\begin{equation}\label{martingale-1}
\mathbb{E} \sup_{0\leq t\leq T}\rho_{t}^{r} \leq
C\Bigl(\mathbb{E}\sup_{0\leq t\leq
T}M_{t}^{p}\Bigr)^{{1/p}}.
\end{equation}
By Doob's maximal inequality, we have%
%
%e2.9 ###
%
\begin{equation}\label{martingale-2}
\EE\sup_{0\leq t\leq T}M_{t}^{p}\leq c_{p}\EE M_{T}^{p}
\end{equation}
for some constant $c_p>0$ depending only on $p$. Finally, choosing
any $\gamma>1$, $\lambda>1$ such that $\frac{1}{\gamma
}+\frac{1}{\lambda}=1$ and applying again the H\"{o}lder inequality
yield
\begin{eqnarray*}
\EE M_{T}^{p} &=&\EE\biggl( \exp\biggl\{ rp\int_{0}^{T}\beta
_{s}\,dW_{s}-\frac{%
\gamma}{2}p^{2}r^{2}\int_{0}^{T}\beta_{s}^{2}\,ds\biggr\} \\
&&\hspace*{56.5pt}{} \times\exp\biggl\{ \frac{\gamma p -1}{2}
pr^{2}\int_{0}^{T}\beta
_{s}^{2}\,ds\biggr\} \biggr) \\
&\leq&\biggl( \EE\exp\biggl\{ rp\gamma\int_{0}^{T}\beta
_{s}\,dW_{s}-\frac{1}{2}%
\gamma^{2}p^{2}r^{2}\int_{0}^{T}\beta_{s}^{2}\,ds\biggr\} \biggr)
^{{1}/{%
\gamma}} \\
&&{}\times\biggl( \EE\exp\biggl\{ \frac{\lambda(\gamma p -1)}{2}
pr^{2}\int_{0}^{T}\beta_{s}^{2}\,ds\biggr\}\biggr) ^{
{1/\lambda}} \\
&=&\biggl( \EE\exp\biggl\{ \frac{\lambda(\gamma p
-1)}{2}pr^{2}\int_{0}^{T}\beta_{s}^{2}\,ds\biggr\} \biggr)
^{{1/\lambda}}<\infty.
\end{eqnarray*}
Combining this inequality with (\ref{martingale-1}) and
(\ref{martingale-2}) we complete the proof.
\end{pf}
\begin{pf*}{Proof of Theorem \ref{T.3.2}} The existence and uniqueness
is well known. We are going to prove
(\ref{e.2.4}).
Let $t\in[0,T]$. Denote $\gamma_{t}=\rho_{t}^{-1}$, where
$\{\rho_t\}_{0\le t\le T}$ is the solution to
(\ref{joint}). Then $\{\gamma_{t}\}_{0\le t\le T}$ satisfies the
following linear stochastic differential equation:
\[
\cases{
d\gamma_{t}=(-\alpha_{t}+\beta_{t}^{2})\gamma_{t}\,dt-\beta
_{t}\gamma
_{t}\,dW_{t}, &\quad$0\le t\le T$,\cr
\gamma_{0}=1.}
\]
For any $0\leq s\leq t\leq T$ and any positive number $r\geq1$,
we have, using~(H2), the H\"{o}lder inequality, the
Burkholder--Davis--Gundy inequality and Lem\-ma~\ref{l.3.1} applied to the
process $\{\gamma_{t}\}_{0\le t\le T}$,
%
%e2.10 ###
%
\begin{eqnarray}\label{e.3.9}
\mathbb{E}|\gamma_{t}-\gamma_{s}|^{r}
&=& \mathbb{E}\biggl| \int_{s}^{t}(-\alpha_{u}+\beta_{u}^{2})\gamma
_{u}\,du-\int_{s}^{t}\beta_{u}\gamma_{u}\,dW_{u}\biggr|^{r} \nonumber\\
&\leq&2^{r-1}\biggl[ {{\mathbb{E}}}\biggl| \int_{s}^{t}(-\alpha
_{u}{+}%
\beta_{u}^{2})\gamma_{u}\,du\biggr| ^{r}{+}C_{r}{{\mathbb{E}}}
\biggl| {{%
\int_{s}^{t}}}\beta_{u}^{2}\gamma_{u}^{2}\,du\biggr| ^{
{r/2}}\biggr]
\\
&\leq&C(t-s)^{r/2},\nonumber
\end{eqnarray}
where $C_r$ is a constant depending only on $r$, and $C$ is a
constant depending on $T$, $r$ and the constants appearing in
conditions (H1) and (H2).

From (\ref{joint}), (\ref{linear}) and by It\^{o}'s formula, we
obtain
\[
d(Y_{t}\rho_{t})=-\rho_{t}f_{t}\,dt+(\beta_{t}\rho_{t}Y_{t}+\rho
_{t}Z_{t})\,dW_{t} .
\]
As a consequence,
%
%e2.11 ###
%
\begin{equation}\label{e.3.10}\qquad
Y_{t} = \rho_{t}^{-1}\mathbb{E}\biggl( \xi\rho_{T}+\int
_{t}^{T}\rho
_{r}f_{r}\,dr\Big|\mathcal{F}_{t}\biggr)
= \mathbb{E}\biggl( \xi\rho_{t,T}+\int_{t}^{T}\rho
_{t,r}f_{r}\,dr\Big|\mathcal{F%
}_{t}\biggr) ,
\end{equation}
where we write $\rho_{t,r}=\rho_{t}^{-1}\rho_{r}=\gamma
_{t}\rho_{r}$ for any $0\le t\leq r\leq T$.

Now, fix $0\leq s\leq t\leq T$. We have
\begin{eqnarray*}
\mathbb{E}|Y_{t}-Y_{s}|^{p}
&=&\mathbb{E} \biggl| \EE\biggl( \xi\rho_{t,T}+\int_{t}^{T}\rho
_{t,r}f_{r}\,dr\Big|\mathcal{F}_{t}\biggr) -\mathbb{E}\biggl( \xi
\rho_{s,T}+\int_{s}^{T}\rho
_{s,r}f_{r}\,dr\Big|\mathcal{F}_{s}\biggr) \biggr|
^{p} \\
&\leq&2^{p-1}\biggl[\mathbb{E}\bigl|\mathbb{E}( \xi\rho
_{t,T}|\mathcal{F}%
_{t}) -\mathbb{E}(\xi\rho_{s,T}|\mathcal{F}_{s})
\bigr|^{p} \\
&&\hspace*{25.4pt}{}+\mathbb{E}\biggl| \mathbb{E}\biggl( \int_{t}^{T}\rho
_{t,r}f_{r}\,dr\Big|\mathcal{%
F}_{t}\biggr) -\mathbb{E}\biggl( \int_{s}^{T}\rho_{s,r}f_{r}\,dr
\Big|\mathcal{F}%
_{s}\biggr) \biggr| ^{p}\biggr] \\
&=&2^{p-1}(I_{1}+I_{2}) .
\end{eqnarray*}
First we estimate $I_{1}$. We have
\begin{eqnarray*}
I_{1} &=&\mathbb{E}\bigl|\mathbb{E}( \xi\rho_{t,T}
|\mathcal{F}_{t}) -%
\mathbb{E}( \xi\rho_{s,T}|\mathcal{F}_{s})
\bigr|^{p} \\
&=&\mathbb{E}\bigl|\mathbb{E}( \xi\rho_{t,T}|\mathcal
{F}_{t}) -%
\mathbb{E}( \xi\rho_{s,T}|\mathcal{F}_{t})
+\mathbb{E}(
\xi\rho_{s,T}|\mathcal{F}_{t}) -\mathbb{E}( \xi
\rho_{s,T}|%
\mathcal{F}_{s}) \bigr|^{p} \\
&\leq&2^{p-1}\bigl[ \mathbb{E}\bigl|\EE( \xi\rho_{t,T}
|\mathcal{F}%
_{t}) -\mathbb{E}( \xi\rho_{s,T}|\mathcal
{F}_{t}) |^{p}+%
\mathbb{E}|\mathbb{E}( \xi\rho_{s,T}|\mathcal
{F}_{t}) -\mathbb{E%
}( \xi\rho_{s,T}|\mathcal{F}_{s})
\bigr|^{p}\bigr] \\
&\leq&2^{p-1}\bigl[ \mathbb{E}|\xi(\rho_{t,T}-\rho
_{s,T})|^{p}+\mathbb{E}%
\bigl|\mathbb{E}( \xi\rho_{s,T}|\mathcal{F}_{t})
-\mathbb{E}(
\xi\rho_{s,T}|\mathcal{F}_{s}) \bigr|^{p}\bigr] \\
&=&2^{p-1}(I_{3}+I_{4}).
\end{eqnarray*}
Using the H\"{o}lder inequality, Lemma \ref{l.3.1} and the estimate
(\ref{e.3.9}) with $r=\frac{%
2pq}{q-p}$, the term $I_{3}$ can be estimated as follows:
\begin{eqnarray*}
I_{3} &\leq&( \mathbb{E}|\xi|^{q}) ^{{p}/{q}}
\bigl( \mathbb{%
E}|\rho_{t,T}-\rho_{s,T}|^{{pq}/({q-p})}\bigr) ^{({q-p})/{q}}
\\
&\leq&( \mathbb{E}|\xi|^{q}) ^{{p}/{q}}\bigl(
\mathbb{E}%
|\gamma_{t}-\gamma_{s}|^{{2pq}/({q-p})}\bigr) ^{
({q-p})/({2q})}\bigl(
\mathbb{E}\rho_{T}^{{2pq}/({q-p})}\bigr) ^{({q-p})/({2q})}\\
&\leq& C|t-s|^{{p}/{2}},
\end{eqnarray*}
where $C$ is a constant depending only on $p,q,T$, $\mathbb{E}|\xi|^{q}$
and the constants appearing in conditions (H1) and (H2).

In order to estimate the term $I_{4}$ we will make use of the
condition $\xi\rho_{T}\in M^{2,q}$. This condition implies that
\[
\xi\rho_{T}=\mathbb{E}(\xi\rho_{T})+\int_{0}^{T}u_{r}\,dW_{r},
\]
where $u$ is a progressively measurable process satisfying \mbox{$\sup
_{0\leq t\leq T}\mathbb{E}%
|u_{t}|^{q}<\infty$}. Therefore, by the Burkholder--Davis--Gundy
inequality, we have
\begin{eqnarray*}
&&
\mathbb{E}\bigl|\mathbb{E}(\xi\rho
_{T}|\mathcal{F}_{t})-\mathbb{E}(\xi\rho
_{T}|\mathcal{F}_{s})\bigr|^{q}\\
&&\qquad=\mathbb{E}\biggl|
\int_{s}^{t}u_{r}\,dW_{r}\biggr| ^{q}
\leq C_{q}\mathbb{E}\biggl| \int_{s}^{t}u_{r}^{2}\,dr\biggr|
^{{q}/{2}%
} \\
&&\qquad\leq C_{q}(t-s)^{({q-2})/{2}}\mathbb{E}\biggl(
\int_{s}^{t}|u_{r}|^{q}\,dr\biggr) \\
&&\qquad\leq C_{q}(t-s)^{{q}/{2}}\sup_{0\leq t\leq T}\mathbb{E}|u_{t}|^{q}.
\end{eqnarray*}
As a consequence, from the definition of $I_{4}$ we have
\begin{eqnarray*}
I_{4} &=&\mathbb{E}\bigl|\gamma_{s}[\mathbb{E}(\xi\rho_{T}|\mathcal
{F}_{t})-%
\mathbb{E}(\xi\rho_{T}|\mathcal{F}_{s})]\bigr|^{p} \\
&\leq&\bigl( \mathbb{E}\gamma_{s}^{{pq}/({q-p})}\bigr) ^{
({q-p})/{q}%
}\bigl( \mathbb{E}\bigl|\mathbb{E}(\xi\rho_{T}|\mathcal
{F}_{t})-\mathbb{E}(\xi
\rho_{T}|\mathcal{F}_{s})\bigr|^{q}\bigr) ^{{p/q}}
\\
&\leq& C|t-s|^{{p/2}},
\end{eqnarray*}
where $C$ is a constant depending on $p,q,T, \sup_{0\leq t\leq T}%
\mathbb{E}|u_{t}|^{q}<\infty$ and the constants appearing in conditions
(H1) and (H2).

The term $I_{2}$ can be decomposed as follows:
\begin{eqnarray*}
I_{2} &=&\mathbb{E}\biggl| \mathbb{E}\biggl( \int_{t}^{T}\rho
_{t,r}f_{r}\,dr\Big|%
\mathcal{F}_{t}\biggr) -\mathbb{E}\biggl( \int_{s}^{T}\rho
_{s,r}f_{r}\,dr\Big|%
\mathcal{F}_{s}\biggr)\biggr| ^{p} \\
&\leq&3^{p-1}\biggl[\mathbb{E}\biggl| \mathbb{E}\biggl( \int
_{t}^{T}\rho
_{t,r}f_{r}\,dr\Big|\mathcal{F}_{t}\biggr) -\mathbb{E}\biggl( \int
_{t}^{T}\rho
_{s,r}f_{r}\,dr\Big|\mathcal{F}_{t}\biggr) \biggr| ^{p} \\
&&\hphantom{3^{p-1}\biggl[}
{} +\mathbb{E}\biggl| \mathbb{E}\biggl( \int_{t}^{T}\rho
_{s,r}f_{r}\,dr\Big|\mathcal{%
F}_{t}\biggr) -\mathbb{E}\biggl( \int_{s}^{T}\rho
_{s,r}f_{r}\,dr\Big|\mathcal{F}%
_{t}\biggr) \biggr| ^{p} \\
&&\hphantom{3^{p-1}\biggl[}{} +\mathbb{E}\biggl| \mathbb{E}\biggl( \int_{s}^{T}\rho
_{s,r}f_{r}\,dr\Big|\mathcal{%
F}_{t}\biggr) -\mathbb{E}\biggl( \int_{s}^{T}\rho
_{s,r}f_{r}\,dr\Big|\mathcal{F}%
_{s}\biggr) \biggr| ^{p}\biggr] \\
&=&3^{p-1}(I_{5}+I_{6}+I_{7}) .
\end{eqnarray*}
Let us first estimate the term $I_{5}$. Suppose that $p<p^{\prime
}<q$. Then, using (\ref{e.3.9}) and the H\"{o}lder inequality, we
can write
\begin{eqnarray*}
I_{5} &=&\mathbb{E}\biggl| \mathbb{E}\biggl( \int_{t}^{T}\rho
_{t,r}f_{r}\,dr\Big|%
\mathcal{F}_{t}\biggr) -\mathbb{E}\biggl( \int_{t}^{T}\rho
_{s,r}f_{r}\,dr\Big|%
\mathcal{F}_{t}\biggr) \biggr| ^{p} \\
&\leq&\mathbb{E}\biggl| \int_{t}^{T}(\rho_{t,r}-\rho
_{s,r})f_{r}\,dr\biggr|
^{p}
= \mathbb{E}\biggl( |\gamma_{t}-\gamma_{s}|^{p}\biggl| \int
_{t}^{T}\rho
_{r}f_{r}\,dr\biggr| ^{p}\biggr) \\
&\leq&\bigl\{ \mathbb{E}| \gamma_{t}-\gamma_{s}|
^{{%
pp^{\prime}}/({p^{\prime}-p})}\bigr\}^{({p^{\prime
}-p})/{p^{\prime}}%
}\biggl\{ \mathbb{E}\biggl| \int_{t}^{T}\rho_{r}f_{r}\,dr\biggr|
^{p^{\prime
}}\biggr\} ^{{p}/{p^{\prime}}} \\
&\leq&C|t-s|^{{p}/{2}}\biggl\{ \mathbb{E}\biggl( \int
_{t}^{T}\rho
_{r}^{2}\,dr\biggr) ^{{p^{\prime}q}/({2(q-p^{\prime})})}\biggr\}
^{{%
p(q-p^{\prime})}/({p^{\prime}q})}\\
&&{}\times\biggl\{ \mathbb{E}\biggl(
\int_{t}^{T}f_{r}^{2}\,dr\biggr) ^{{q}/{2}}\biggr\}^{
{p}/{q}} \\[-2pt]
&\leq&\widehat{C}|t-s|^{{p}/{2}}\Vert f\Vert_{H^{q}}^{p},
\end{eqnarray*}
where $\widehat{C}$ is a constant depending on $p,p^{\prime}$, $q$, $T$
and the constants appearing in conditions (H1) and (H2).

Now we estimate $I_{6}$. Suppose that $p<p^{\prime}<q$. We have, as in the
estimate of the term $I_{5}$,\vspace*{-2pt}
\begin{eqnarray*}
I_{6} &=&\mathbb{E}\biggl| \mathbb{E}\biggl( \int_{t}^{T}\rho
_{s,r}f_{r}\,dr\Big|%
\mathcal{F}_{t}\biggr) -\mathbb{E}\biggl( \int_{s}^{T}\rho
_{s,r}f_{r}\,dr\Big|%
\mathcal{F}_{t}\biggr) \biggr| ^{p} \\[-2pt]
&\leq&\mathbb{E}\biggl| \int_{s}^{t}\rho_{s,r}f_{r}\,dr\biggr|
^{p}=\mathbb{E}%
\biggl( \rho_{s}^{-p}\biggl| \int_{s}^{t}\rho_{r}f_{r}\,dr\biggr|
^{p}\biggr)
\\[-2pt]
&\leq&\bigl\{ \mathbb{E}\rho_{s}^{-{pp^{\prime}}/({p^{\prime
}-p})%
}\bigr\}^{({p^{\prime}-p})/{p^{\prime}}}\biggl\{ \mathbb
{E}\biggl|
\int_{s}^{t}\rho_{r}f_{r}\,dr\biggr| ^{p^{\prime}}\biggr\} ^{
{p}/{%
p^{\prime}}} \\[-2pt]
&=&C\biggl\{ \mathbb{E}\biggl| \int_{s}^{t}\rho_{r}f_{r}\,dr\biggr|
^{p^{\prime
}}\biggr\} ^{{p}/{p^{\prime}}} \\[-2pt]
&\leq&C|t-s|^{{p}/{2}}\Bigl\{ \mathbb{E}\sup_{0\leq t\leq
T}\rho_{t}^{%
{p^{\prime}q}/({q-p^{\prime}})}\Bigr\} ^{{p(q-p^{\prime
})}/({%
p^{\prime}q})}\Vert f\Vert_{H^{q}}^{p}\\[-2pt]
&=& \widehat{C}|t-s|^{{p}/{2}} ,\vspace*{-2pt}
\end{eqnarray*}
where $\widehat{C}$ is a constant depending on $p,p^{\prime}$,
$q$, $T$ and the constants appearing in conditions (H1) and (H2).

The fact that $\int_{0}^{T}\rho
_{r}f_{r}\,dr$ belongs to $M^{2,q}$ implies that
\[
\int_{0}^{T}\rho_{r}f_{r}\,dr=\mathbb{E}\int_{0}^{T}\rho
_{r}f_{r}\,dr+\int_{0}^{T}v_{r}\,dW_{r},\vspace*{-2pt}
\]
where $\{v_t\}_{0\le t\le T}$ is a progressively measurable process
satisfying
\[
\sup_{0\leq t\leq T}\mathbb{E}%
|v_{t}|^{q}<\infty.\vspace*{-2pt}
\]
Then, by the Burkholder--Davis--Gundy
inequality we have\vspace*{-2pt}
\begin{eqnarray*}
&&\mathbb{E}\biggl| \mathbb{E}\biggl( \int_{s}^{T}\rho
_{r}f_{r}\,dr\Big|\mathcal{F}%
_{t}\biggr) -\mathbb{E}\biggl( \int_{s}^{T}\rho_{r}f_{r}\,dr\Big|\mathcal
{F}%
_{s}\biggr) \biggr| ^{q}\\[-2pt]
&&\qquad=%
\EE\biggl|\mathbb{E}\biggl( \int_{0}^{T}\rho_{r}f_{r}\,dr\Big|\mathcal
{F}_{t}\biggr) -%
\mathbb{E}\biggl( \int_{0}^{T}\rho_{r}f_{r}\,dr\Big|\mathcal{F}_{s}
\biggr)\biggr|^q \\[-2pt]
&&\qquad=\mathbb{E}\biggl| \int_{s}^{t}v_{r}\,dW_{r}\biggr| ^{q} \leq
C_{q}(t-s)^{{q}/{2}}\sup_{0\leq t\leq T}\mathbb{E}|v_{t}|^{q}.\vspace*{-2pt}
\end{eqnarray*}
Finally, we estimate $I_{7}$ as follows:
%
%e2.12 ###
%
\begin{eqnarray}\label{e.3.9a}
I_{7} &=&\mathbb{E}\biggl| \mathbb{E}\biggl( \int_{s}^{T}\rho
_{s,r}f_{r}\,dr\Big|%
\mathcal{F}_{t}\biggr) -\mathbb{E}\biggl( \int_{s}^{T}\rho
_{s,r}f_{r}\,dr\Big|%
\mathcal{F}_{s}\biggr) \biggr| ^{p} \nonumber\\
&=&\mathbb{E}\biggl| \rho_{s}^{-1}\biggl( \mathbb{E}\biggl( \int
_{s}^{T}\rho
_{r}f_{r}\,dr\Big|\mathcal{F}_{t}\biggr) -\mathbb{E}\biggl( \int
_{s}^{T}\rho
_{r}f_{r}\,dr\Big|\mathcal{F}_{s}\biggr) \biggr) \biggr| ^{p} \nonumber\\
&\leq&\bigl\{
\mathbb{E}\rho_{s}^{-{pq}/({q-p})}\bigr\}
^{({q-p})/{%
p}}\nonumber\\[-8pt]\\[-8pt]
&&{}\times\biggl\{ \mathbb{E}\biggl| \mathbb{E}\biggl( \int_{s}^{T}\rho
_{r}f_{r}\,dr\Big|\mathcal{F}_{t}\biggr) -\mathbb{E}\biggl(
\int_{s}^{T}\rho_{r}f_{r}\,dr\Big|\mathcal{F}_{s}\biggr) \biggr|
^{q}\biggr\}^{{p}/{q}}
\nonumber\\
&\leq&C\biggl\{ \mathbb{E}\biggl| \mathbb{E}\biggl( \int
_{s}^{T}\rho
_{r}f_{r}\,dr\Big|\mathcal{F}_{t}\biggr) -\mathbb{E}\biggl(
\int_{s}^{T}\rho_{r}f_{r}\,dr\Big|\mathcal{F}_{s}\biggr) \biggr|
^{q}\biggr\}^{{p}/{q}}\nonumber\\
&\le&\widehat{C}|t-s|^{{p}/{2}} ,\nonumber
\end{eqnarray}
where $\widehat{C}$ is a constant depending on $p$, $q$, $T$,
$\sup_{0\leq t\leq T}\mathbb{E}|v_{t}|^{q}$ and the constants
appearing in conditions (H1) and (H2).

As a consequence, we obtain for all $s, t\in[0,T]$
\[
\mathbb{E}|Y_{t}-Y_{s}|^{p}\leq K|t-s|^{{p}/{2}},
\]
where $K$ is a constant independent of $s$ and $t$.
\end{pf*}

%s2.3 ###
\subsection{The Malliavin calculus for BSDEs}
We return to the study of (\ref{bsde}).
The main assumptions we make on the terminal value $\xi$
and generator $f$ are the following:
\begin{assumption}
\label{a.3.2} Fix $2\leq p<\frac{q}{2}$.
\begin{longlist}
\item$\xi\in\mathbb{D}^{2,q}$, and there exists $L>0$,
such that for all $\theta, \theta^{\prime}\in\lbrack0,T]$,
%
%e2.14 ###
%e2.13 ###
%
\begin{eqnarray}
\label{e5}
\mathbb{E}|D_{\theta}\xi-D_{\theta^{\prime}}\xi|^{p}&\leq&
L|\theta-\theta^{\prime}|^{{p}/{2}},
\\
\label{e2}
\sup_{0\leq\theta\leq T}\mathbb{E}|D_{\theta}\xi|^{q}&<&\infty
\end{eqnarray}
and
%
%e2.15 ###
%
\begin{equation}\label{e2-2}
\sup_{0\leq\theta\leq T}\sup_{0\leq u\leq
T}\mathbb{E}|D_{u}D_{\theta}\xi|^{q}<\infty.
\end{equation}

\item The generator $f(t,y,z)$ has continuous and uniformly
bounded first-
and second-order partial derivatives with respect to $y$ and $z$,
and $f(\cdot,0,0)\in H_{ \mathcal{F}}^{q}([0,T])$.

\item Assume that $\xi$ and $f$ satisfy the above conditions
(i) and (ii). Let $(Y,Z)$ be the unique solution to
(\ref{bsde}) with terminal value $\xi$ and generator~$f$. For each
$(y,z)\in\mathbb{R}\times\mathbb{R}$, $f(\cdot,y,z)$, $\partial
_{y}f(\cdot,y,z)$
and $\partial_{z}f(\cdot,y,z)$\vadjust{\goodbreak} belong to~%
$\mathbb{L}_{a}^{1,q}$, and the Malliavin derivatives $Df(\cdot
,y,z)$, $D\partial_{y}f(\cdot,y,z)$ and $D\partial_{z}f(\cdot
,y,z)$ satisfy
%
%e2.18 ###
%e2.17 ###
%e2.16 ###
%
\begin{eqnarray}
\label{e3}
\sup_{0\leq\theta\leq T}\EE\biggl( \int_{\theta}^{T} |D_{\theta
}f(t,Y_t,Z_t)|^{2}\,dt\biggr) ^{{q}/{2}}&<&\infty,
\\
\label{e5-1}
\sup_{0\leq\theta\leq T}\mathbb{E}\biggl(
\int_{\theta}^{T}|D_{\theta}\partial
_{y}f(t,Y_t,Z_t)|^{2}\,dt\biggr) ^{{q}/{2}}&<&\infty,
\\
\label{e6}
\sup_{0\leq\theta\leq T}\mathbb{E}\biggl(
\int_{\theta}^{T}|D_{\theta}\partial
_{z}f(t,Y_t,Z_t)|^{2}\,dt\biggr) ^{{q}/{2}}&<&\infty,
\end{eqnarray}
and there exists $L\,{>}\,0$ such that for any $t\,{\in}\,(0,T]$, and for any
\mbox{$0\,{\le}\,\theta, \theta^{\prime}\,{\leq}\, t\,{\le}\, T$}
%
%e2.19 ###
%
\begin{equation} \label{e4}\qquad
\mathbb{E}\biggl( \int_{t}^{T}|D_{\theta}f(r,Y_r,Z_r)-D_{\theta
^{\prime}}f(r,Y_r,Z_r)|^{2}\,dr\biggr) ^{{p}/{2}}\leq L|\theta
-\theta^{\prime}|^{{p}/{2}}.
\end{equation}
For each $\theta\in\lbrack0,T]$, and each pair of $(y,z)$,
$D_{\theta}f(\cdot,y,z)\in\mathbb{L}_{a}^{1,q}$ and it has
continuous partial derivatives
with respect to $y,z$, which are denoted by $\partial_yD_\theta
f(t,y,z)$
and $\partial_zD_\theta f(t,y,z)$, and the Malliavin derivative
$D_{u}D_{\theta}f(t,y,z)$ satisfies
%
%e2.20 ###
%
\begin{equation}\label{e6-2}
\sup_{0\leq\theta\leq T}\sup_{0\leq u\leq T} \mathbb{E}\biggl(
\int_{\theta\vee u}^{T}|D_{u}D_{\theta}f(t,Y_t,Z_t)|^{2}\,dt\biggr)
^{{q}/{2}}<\infty.
\end{equation}
\end{longlist}
\end{assumption}

The following property is easy to check and we omit the proof.
\begin{remark} Conditions (\ref{e5-1}) and (\ref{e6}) imply
\[
\sup_{0\leq\theta\leq T}\mathbb{E}\biggl(
\int_{\theta}^{T}|\partial_{y}D_{\theta
}f(t,Y_t,Z_t)|^{2}\,dt\biggr) ^{{q}/{2}}<\infty
\]
and
\[
\sup_{0\leq\theta\leq T}\mathbb{E}\biggl(
\int_{\theta}^{T}|\partial_{z}D_{\theta
}f(t,Y_t,Z_t)|^{2}\,dt\biggr) ^{{q}/{2}}<\infty,
\]
respectively.
\end{remark}

The following is the main result of this section.
\begin{theorem}\label{t.3.1}
Let Assumption \ref{a.3.2} be satisfied.
\begin{longlist}[(a)]
\item[(a)] There exists a unique solution pair $\{(Y_t,Z_t)\}_{
0\le t\le T}$ to the BSDE (\ref{bsde}), and $Y, Z$ are in
$ \mathbb{L}^{1,q}_a$. A version of the Malliavin derivatives
$\{(D_\theta Y_t,\break D_\theta Z_t)\}_{ 0\leq\theta, t\leq T}$ of the
solution pair
satisfies the following linear BSDE:
%
%e2.22 ###
%e2.21 ###
%
\begin{eqnarray}
\label{e.3.12}\qquad
D_\theta Y_t &=&D_\theta\xi+\int_t^T[\partial_y{f(r,Y_r,Z_r)}
D_\theta
Y_r \nonumber\\[-2pt]
&&\hphantom{D_\theta\xi+\int_t^T[}
{}+\partial_z{f(r,Y_r,Z_r)}D_\theta Z_r +D_\theta f(r,Y_r,Z_r)]\,dr\\
&&{} -\int_t^TD_\theta
Z_r\,dW_r,\qquad 0\le\theta\leq t\leq T ;\nonumber\vadjust{\goodbreak}\\[-2pt]
\label{e.3.12-2}
D_\theta Y_t&=&0,\qquad D_\theta Z_t=0,\qquad 0\le t<\theta\leq
T.
\end{eqnarray}
Moreover, $\{D_{t} Y_t\}_{0\leq t\leq T}$ defined by
(\ref{e.3.12}) gives a version of $\{Z_t\}_{0\leq t\leq T}$,
namely, $\mu\times P $ a.e.
%
%e2.23 ###
%
\begin{equation}\label{e.3.13}
Z_t=D_{t} Y_t .
\end{equation}
\item[(b)] There exists a constant $K>0$, such that, for all $%
s, t\in[0,T]$,
%
%e2.24 ###
%
\begin{equation}\label{e-z}
\mathbb{E}|Z_{t}-Z_{s}|^{p}\leq K|t-s|^{{p}/{2}}.
\end{equation}
\end{longlist}
\end{theorem}
\begin{pf} Part (a): The proof of the existence and uniqueness of the
solution $(Y,Z)$, and $Y, Z\in
\mathbb{L}_{a}^{1,2}$ is similar to that of Proposition 5.3 in~%
\cite{KPQ97}, and also the fact that $(D_{\theta
}Y_{t},D_{\theta}Z_{t})$ is given by (\ref{e.3.12}) and
(\ref{e.3.12-2}). In Proposition~5.3 in~\cite{KPQ97} the
exponent $q$ is equal to $4$, and one assumes that
$\int_{0}^{T}\Vert D_{\theta}f(\cdot,Y,\allowbreak Z)\Vert
_{H^{2}}^{2}\,d\theta<\infty$, which
is a consequence of (\ref{e3}) and the fact that $Y, Z\in
\mathbb{L}%
_{a}^{1,2}$.

Furthermore, from conditions (\ref{e2}) and (\ref{e3}) and the
estimate in Lem\-ma~\ref{T.2.1}, we obtain
%
%e2.25 ###
%
\begin{equation}\label{e.3.19-1}
\sup_{0\le\theta\le T}\biggl\{\EE\sup_{\theta\le t\le T}|D_\theta
Y_t|^q+\EE\biggl(\int_\theta^T|D_\theta
Z_t|^2\,dt\biggr)^{{q}/{2}}\biggr\}<\infty.
\end{equation}
Hence, by Proposition 1.5.5 in \cite{N06}, $Y$ and $Z$ belong to
$\mathbb{L}_a^{1,q}$.

Part (b): Let $0\le s\leq t\leq T$. In this proof, $C>0$ will be a
constant independent of $s$ and $t$, and may vary from line to
line.

By representation (\ref{e.3.13}) we have
%
%e2.26 ###
%
\begin{equation}\label{e.3.18}
Z_t-Z_s=D_tY_t-D_sY_s=(D_tY_t-D_sY_t)+(D_sY_t-D_sY_s).
\end{equation}
From Lemma \ref{T.2.1} and equation (\ref{e.3.12}) for $\theta=s$
and $\theta'=t$, respectively, we obtain, using conditions
(\ref{e5}) and (\ref{e4}),
%
%e2.27 ###
%
\begin{eqnarray}\label{e.3.19}
&&\EE|D_tY_t-D_sY_t|^p+\EE\biggl(\int_t^T|D_tZ_r-D_sZ_r|^2\,dr
\biggr)^{{p}/{2}}\nonumber\\[-2pt]
&&\qquad\le
C\biggl[\EE|D_t\xi-D_s\xi|^p\nonumber\\[-10pt]\\[-10pt]
&&\qquad\quad\hphantom{C\biggl[}
{}+\EE\biggl(\int
_t^T|D_tf(r,Y_r,Z_r)-D_sf(r,Y_r,Z_r)|^2\,dr\biggr)^{{p}/{2}}
\biggr]\nonumber\\[-2pt]
&&\qquad\le C|t-s|^{{p}/{2}}.\nonumber
\end{eqnarray}
Denote $\alpha_u=\partial_y f(u,Y_u,Z_u)$ and $\beta_u=\partial_z
f(u,Y_u,Z_u)$ for all $u\in[0,T]$. Then, by Assumption \ref{a.3.2}(ii),
the processes $\alpha$ and $\beta$ satisfy conditions (H1)
and~(H2) in Assumption \ref{a.2.1}, and from (\ref{e.3.12}) we
have for $r\in[s,T]$
\[
D_sY_r=D_s\xi+\int_r^T[\alpha_uD_sY_u+\beta
_uD_sZ_u+D_sf(u,Y_u,Z_u)]\,du-\int_r^TD_sZ_u\,dW_u.
\]
Next, we are going to use Theorem \ref{T.3.2} to estimate
$\EE|D_sY_t-D_sY_s|^p$. Fix $p^\prime$ with
$p<p^\prime<\frac{q}{2}$ (notice that $p^\prime<\frac{q}{2}$ is
equivalent to $\frac{p^\prime}{q-p^\prime}<1$). From conditions~
(\ref{e2}) and (\ref{e3}), it is obvious that $D_s\xi\in
L^q(\Omega)\subset L^{p^\prime}(\Omega)$ and $D_s f(\cdot,Y,Z)\in
H^q([0,T])\subset H^{p^\prime}([0,T])$ for any $s\in[0,T]$.
We are going to show
that, for any $s\in[0,T]$, $\rho_TD_s\xi$ and $\int_s^T \rho_u
D_sf(u,Y_u,Z_u)\,du$ are
elements in $M^{2,p^\prime}$, where
\[
\rho_r=\exp\biggl\{\int_0^r\beta_u\,dW_u+\int_0^r\biggl(\alpha
_u-\frac{1}{2}\beta_u^2\biggr)\,du\biggr\}.
\]
For any $0\le\theta\le r\le T$, let us compute
\begin{eqnarray*}
D_\theta\rho_r&=&\rho_r\biggl\{\int_\theta^r[\partial
_{yz}f(u,Y_u,Z_u)D_\theta
Y_u\\[-2pt]
&&\hphantom{\rho_r\biggl\{\int_\theta^r[}
{} +\partial_{zz}f(u,Y_u,Z_u)D_\theta
Z_u+ D_\theta\partial_{z}f(u,Y_u,Z_u)]\,dW_u\\[-2pt]
&&\hphantom{\rho_r\biggl\{}
{}+\partial_zf(\theta
,Y_\theta,Z_\theta)\\[-2pt]
&&\hphantom{\rho_r\biggl\{}+\int_\theta^r\bigl(\partial_{yy}f(u,Y_u,Z_u)-\partial
_{yz}f(u,Y_u,Z_u)\beta_u\bigr)D_\theta
Y_u\,du\\[-2pt]
&&\hphantom{\rho_r\biggl\{}
{}+\int_\theta^r\bigl(\partial_{yz}f(u,Y_u,Z_u)-\partial
_{zz}f(u,Y_u,Z_u)\beta_u\bigr)D_\theta
Z_u\,du\\[-2pt]
&&\hphantom{\rho_r\biggl\{}
\hspace*{17.3pt}{}+\int_\theta^r\bigl(D_\theta\partial_{y}f(u,Y_u,Z_u)-\beta_uD_\theta
\partial_{z}f(u,Y_u,Z_u)\bigr)\,du\biggr\}.
\end{eqnarray*}
By the boundedness of the first- and second-order partial
derivatives of~$f$ with respect to $y$ and $z$, (\ref{e5-1}),
(\ref{e6}), (\ref{e.3.19-1}), Lemma \ref{l.3.1}, the H\"{o}lder
inequality and the Burkholder--Davis--Gundy inequality, it is easy to
show that for any $p''<q$,
%
%e2.28 ###
%
\begin{equation}\label{e.3.21}
\sup_{0\le\theta\le T}\EE\sup_{\theta\le r\le
T}|D_\theta\rho_r|^{p''}<\infty.
\end{equation}
By the Clark--Ocone--Haussman formula, we have
\begin{eqnarray*}
\rho_TD_s\xi&=&\EE(\rho_TD_s\xi)+\int_0^T\EE(D_\theta(\rho
_TD_s\xi)|\mathcal{F}_\theta)
\,dW_\theta\\[-2pt]
&=&\EE(\rho_TD_s\xi)+\int_0^T\EE(D_\theta\rho_TD_s\xi+\rho
_TD_\theta
D_s\xi|\mathcal{F}_\theta) \,dW_\theta\\[-2pt]
&=&\EE(\rho_TD_s\xi)+\int_0^Tu_\theta^s \,dW_\theta
\end{eqnarray*}
and
\begin{eqnarray*}
&&\int_s^T\rho_rD_sf(r,Y_r,Z_r)\,dr\\
&&\qquad=\EE\int_s^T\rho_rD_sf(r,Y_r,Z_r)\,dr\\
&&\qquad\quad{}+\int_0^T\EE\biggl(D_\theta
\int_s^T\rho_rD_sf(r,Y_r,Z_r)\,dr\Big|\mathcal{F}_\theta
\biggr)\,dW_\theta\\
&&\qquad=\EE\int_s^T\rho_rD_sf(r,Y_r,Z_r)\,dr\\
&&\qquad\quad{}+\int_0^T\EE\biggl(\int_s^T[D_\theta\rho_rD_sf(r,Y_r,Z_r)\\
&&\hphantom{{}+\int_0^T\EE\biggl(\int_s^T[}
\qquad\quad{}+\rho
_r\partial_yD_sf(r,Y_r,Z_r)D_\theta
Y_r\\
&&\hphantom{{}+\int_0^T\EE\biggl(\int_s^T[}
\qquad\quad{}+\rho_r\partial_zD_sf(r,Y_r,Z_r)D_\theta Z_r\\
&&\qquad\quad\hspace*{86.12pt}{}+
\rho_rD_\theta
D_sf(r,Y_r,Z_r)]\,dr\Big|\mathcal{F}_\theta\biggr)\,dW_\theta\\
&&\qquad=\EE\int_s^T\rho_rD_sf(r,Y_r,Z_r)\,dr+\int_0^Tv_\theta^s
\,dW_\theta.
\end{eqnarray*}
We claim that $\sup_{0\le\theta\le
T}\EE|u^s_\theta|^{p^{\prime}}<\infty$ and $\sup_{0\le\theta\le
T}\EE|v^s_\theta|^{p^{\prime}}<\infty$. In fact,
\begin{eqnarray*}
\EE|u^s_\theta|^{p^{\prime}}
&=&\EE\bigl|\EE(D_\theta\rho_TD_s\xi+\rho_TD_\theta
D_s\xi|\mathcal{F}_\theta)\bigr|^{p^{\prime}}\\
&\le&2^{{p^{\prime}}-1}(\EE|D_\theta\rho_TD_s\xi
|^{p^{\prime}}+\EE|\rho_TD_\theta
D_s\xi|^{p^{\prime}})\\
&\le&2^{{p^{\prime}}-1}\bigl(\bigl(\EE|D_\theta\rho_T|^{
{p^{\prime}q}/({q-p^{\prime}})}\bigr)^{({q-p^{\prime}})/{q}}%
(\EE|D_s\xi|^q)^{{p^{\prime}}/{q}}\\
&&\hspace*{27.4pt}{}+\bigl(\EE\rho_T^{{p^{\prime}q}/({q-p^{\prime}})}
\bigr)^{({q-p^{\prime}})/{q}}%
(\EE|D_\theta D_s\xi|^q)^{{p^{\prime}}/{q}}\bigr).
\end{eqnarray*}
By (\ref{e2}), (\ref{e2-2}), (\ref{e.3.21}) and Lemma \ref{l.3.1},
we have $\sup_{0\le s\le T}\sup_{0\le\theta\le
T}\EE|u^s_\theta|^{p^{\prime}}<\infty$. On the other hand,
\begin{eqnarray*}
\EE|v^s_\theta|^{p^{\prime}}
&=&\EE\biggl\vert\EE\biggl(\int_s^T[D_\theta\rho
_rD_sf(r,Y_r,Z_r)\\
&&\hphantom{\EE\biggl\vert\EE\biggl(\int_s^T[}
{}+\rho_r\partial_yD_sf(r,Y_r,Z_r)D_\theta
Y_r\\
&&\hphantom{\EE\biggl\vert\EE\biggl(\int_s^T[}
{} +\rho_r\partial_zD_sf(r,Y_r,Z_r)D_\theta Z_r\\
&&\hspace*{65.4pt}{}
+\rho_rD_\theta
D_sf(r,Y_r,Z_r)]\,dr\Big|\mathcal{F}_\theta\biggr)\biggr\vert
^{p^{\prime}}\\
&\le&4^{p^{\prime}-1}[J_1+J_2+J_3+J_4],
\end{eqnarray*}
where
\begin{eqnarray*}
J_1&=&\EE\biggl\vert\int_s^TD_\theta\rho_rD_sf(r,Y_r,Z_r)\,dr
\biggr\vert^{p^{\prime}},
\\
J_2&=&\EE\biggl\vert\int_s^T\rho_r\partial_yD_sf(r,Y_r,Z_r)D_\theta
Y_r\,dr\biggr\vert^{p^{\prime}},
\\
J_3&=&\EE\biggl\vert\int_s^T\rho_r\partial_zD_sf(r,Y_r,Z_r)D_\theta
Z_r\,dr\biggr\vert^{p^{\prime}}
\end{eqnarray*}
and
\[
J_4=\EE\biggl\vert\int_s^T\rho_rD_\theta
D_sf(r,Y_r,Z_r)\,dr\biggr\vert^{p^{\prime}}.
\]
For $J_1$, we have
\begin{eqnarray*}
J_1 &\le&\EE\biggl({\sup_{\theta\le r\le T}}\vert
D_\theta\rho_r\vert^{p^{\prime}}\biggl\vert\int
_s^TD_sf(r,Y_r,Z_r)\,dr\biggr\vert^{p^{\prime}}\biggr)\\
&\le&\Bigl({\EE\sup_{\theta\le r\le T}}\vert
D_\theta\rho_r\vert^{{p^{\prime}q}/({q-p^{\prime}})}
\Bigr)^{({q-p^{\prime}})/{q}}%
\\
&&{}\times
\biggl(\EE\biggl\vert\int_s^TD_sf(r,Y_r,Z_r)\,dr\biggr\vert^q
\biggr)^{{p^{\prime}}/{q}}\\
&\le&T^{{p^{\prime}}/{2}}\Bigl({\EE\sup_{\theta\le r\le
T}}\vert
D_\theta\rho_r\vert^{{p^{\prime}q}/({q-p^{\prime}})}
\Bigr)^{({q-p^{\prime}})/{q}}%
\\
&&{}\times\biggl(\EE\biggl(\int_0^T\vert
D_sf(r,Y_r,Z_r)\vert^2\,dr\biggr)^{{q}/{2}}\biggr)^{
{p^{\prime}}/{q}}.
\end{eqnarray*}
For $J_2$, we have
\begin{eqnarray*}
J_2 &\le&\EE\biggl({\sup_{\theta\le r\le T}}\vert D_\theta
Y_r\vert^{p^{\prime}}%
\biggl(\sup_{0\le r\le T}\rho_r%
\int_s^T\vert\partial_yD_sf(r,Y_r,Z_r)\vert
\,dr\biggr)^{p^{\prime}}\biggr)\\%
&\le&\Bigl(\EE\sup_{\theta\le r\le T}\vert D_\theta
Y_r\vert^{q}\Bigr)^{{p^{\prime}}/{q}}\\
&&{}\times
\biggl(\EE\biggl(\sup_{0\le r\le T}\rho_r
\int_s^T\vert\partial_yD_sf(r,Y_r,Z_r)\vert
\,dr\biggr)^{{p^{\prime}q}/({q-p^{\prime}})}\biggr)^{
({q-p^{\prime}})/{q}}\\
&\le&\Bigl({\EE\sup_{\theta\le r\le T}}\vert D_\theta
Y_r\vert^{q}\Bigr)^{{p^{\prime}}/{q}} \Bigl(\EE\sup_{0\le
r\le
T}\rho_r^{{p^{\prime}q}/({q-2p^{\prime}})}\Bigr)^{
({q-2p^{\prime}})/{q}}\\
&&{}\times\biggl(\EE\biggl(\int_s^T\vert\partial
_yD_sf(r,Y_r,Z_r)\vert
\,dr\biggr)^q\biggr)^{{p^{\prime}}/{q}}\\
&\le&T^{{p^{\prime}}/{2}}\Bigl(\EE\sup_{\theta\le r\le
T}\vert
D_\theta
Y_r\vert^{q}\Bigr)^{{p^{\prime}}/{q}}%
\Bigl(\EE\sup_{0\le r\le
T}\rho_r^{{p^{\prime}q}/({q-2p^{\prime}})}\Bigr)^{
({q-2p^{\prime}})/{q}}\\
&&{}\times\biggl(\EE\biggl(\int_0^T\vert\partial
_yD_sf(r,Y_r,Z_r)\vert^2
\,dr\biggr)^{{q}/{2}}\biggr)^{{p^{\prime}}/{q}}.
\end{eqnarray*}
Using a similar techniques as before, we obtain that
\begin{eqnarray*}
J_3&\le&T^{{p'/2}}\biggl(\EE\biggl(\int_0^T\vert D_\theta
Z_r\vert^2 \,dr\biggr)^{{q/2}}\biggr)^{{p^{\prime}/q}}
\Bigl(\EE\sup_{0\le r\le
T}\rho_r^{{p^{\prime}q}/({q-2p^{\prime}})}\Bigr)^{
({q-2p^{\prime}})/{q}}\\
&&{}\times\biggl(\EE\biggl(\int_0^T\vert\partial
_zD_sf(r,Y_r,Z_r)\vert^2
\,dr\biggr)^{{q/2}}\biggr)^{{p^{\prime}/q}}
\end{eqnarray*}
and
\begin{eqnarray*}
J_4&\le&T^{{p^{\prime}}/{2}}\Bigl(\EE\sup_{0\le r\le
T}\rho_r^{{p^{\prime}q}/({q-p^{\prime}})}\Bigr)^{
({q-p^{\prime}})/{q}}\\
&&{}\times\biggl(\EE\biggl(\int_0^T\vert
D_\theta D_sf(r,Y_r,Z_r)\vert^2
\,dr\biggr)^{{q/2}}\biggr)^{{p^{\prime}/q}}.
\end{eqnarray*}
By (\ref{e3}), (\ref{e5-1})--(\ref{e6-2}), (\ref{e.3.21}) and Lemma
\ref{l.3.1}, we obtain that
\[
\sup_{0\le s\le T}\sup_{0\le\theta\le
T}\EE|v^s_\theta|^{p^{\prime}}<\infty.
\]
Therefore, $\rho_T\xi$ and
$\int_0^T\rho_uD_sf(u,Y_u,Z_u)\,du$ belong to $M^{2,p'}$.

Thus by Theorem \ref{T.3.2} with $p<p^\prime$, there is a constant
$C(s)>0$, such that
\[
\EE\vert D_sY_t-D_sY_s\vert^p\le C(s)|t-s|^{{p/2}}
\]
for all $t\in[s, T]$. Furthermore, taking into account the proof
of the estimates $I_k$ ($k=3, 4,\ldots, 7$) in the proof of Theorem
\ref{T.3.2},
we can show that $\sup_{0\le s\le T}C(s)=:C<\infty$. Thus we have
%
%e2.29 ###
%
\begin{equation}\label{e.3.23}
\EE\vert D_sY_t-D_sY_s\vert^p\le C|t-s|^{{p/2}}
\end{equation}
for all $s, t\in[0,T]$. Combining (\ref{e.3.23}) with (\ref
{e.3.18}) and
(\ref{e.3.19}), we obtain that there is a constant $K>0$
independent of $s$ and $t$, such that
\[
\EE\vert Z_t-Z_s\vert^p\le K\vert t-s\vert^{{p/2}}
\]
for all $s, t \in[0,T]$.
\end{pf}
\begin{corollary}\label{l.3.7}
Under the assumptions in Theorem \ref{T.2.1}, let $(Y,Z)\in
S^q_\mathcal{F}([0,T])\times H_\mathcal{F}^q([0,T])$ be the unique
solution pair to (\ref{bsde}). If\break $\sup_{0\le t\le
T}\mathbb{E}|Z_t|^{q}<\infty$, then there exists a constant $C$,
such that, for any $s, t\in[0,T]$,
%
%e2.30 ###
%
\begin{equation}\label{Hy}
\mathbb{E}|Y_t-Y_s|^q\le C|t-s|^{{q/2}}.
\end{equation}
\end{corollary}
\begin{pf}
Without loss of generality we assume $0\le s\le t\le T$. $C>0$ is
a~constant independent of $s$ and $t$, which may vary from line to
line. Since
\[
Y_s=Y_t+\int_s^tf(r,Y_r,Z_r)\,dr-\int_s^tZ_r\,dW_r,
\]
we have, by the Lipschitz condition on $f$,
\begin{eqnarray*}
\mathbb{E}|Y_t-Y_s|^q
&=&\mathbb{E}\biggl|\int_s^tf(r,Y_r,Z_r)\,dr-%
\int_s^tZ_r\,dW_r\biggr|^q \\
&\le&2^{q-1}\biggl(\mathbb{E}\biggl|\int_s^tf(r,Y_r,Z_r)\,dr
\biggr|^q+\mathbb{E}%
\biggl|\int_s^tZ_r\,dW_r\biggr|^q\biggr) \\
&\le&C_q\biggl(|t-s|^{{q/2}}\mathbb{E}\biggl(\int
_s^t|f(r,Y_r,Z_r)|^2\,dr%
\biggr)^{{q/2}}+\mathbb{E}\biggl(\int_s^t|Z_r|^2\,dr
\biggr)^{{q/2}}\biggr) \\
&\le&C\biggl\{|t-s|^{{q/2}}\biggl[\mathbb{E}\biggl(\int
_s^t|Y_r|^2\,dr%
\biggr)^{{q/2}}+\mathbb{E}\biggl(\int_s^t|Z_r|^2\,dr
\biggr)^{{q/2}}\\
&&\hphantom{C\biggl\{|t-s|^{{q/2}}\biggl[}
\hspace*{54.4pt}{} +\mathbb{E}\biggl(\int_s^t|f(r,0,0)|^2\,dr\biggr)^{
{q}/{2}}\biggr]\\
&&\hspace*{133pt}{} + |t-s|^{{q/2}}\sup_{0\le r\le T}\mathbb{E}|Z_r|^{q}\biggr\} \\
&\le&C|t-s|^{{q/2}}.
\end{eqnarray*}
The proof is complete.
\end{pf}
\begin{remark}\label{r.3.8}
From Theorem \ref{t.3.1} we know that $\{(D_\theta Y_t, D_\theta
Z_t)\}_{0\le\theta\le t\le T}$ satisfies equation (\ref{e.3.12})
and $Z_t=D_tY_t$, $\mu\times P$ a.e. Moreover, since
(\ref{e2}) and (\ref{e3}) hold, we can apply the estimate
(\ref{e.2.1}) in Lemma \ref{T.2.1} to the linear BSDE
(\ref{e.3.12}) and deduce $\sup_{0\le t\le
T}\mathbb{E}|Z_t|^{q}<\infty$. Therefore, by Lemma \ref{l.3.7},
the process $Y$ satisfies the inequality (\ref{Hy}). By
Kolmogorov's continuity criterion this implies that $Y$ has
H\"{o}lder continuous trajectories of order $\gamma$ for any
$\gamma<\frac12-\frac1q$.
\end{remark}

%s2.4 ###
\subsection{Examples}
In this section we discuss three particular examples where
Assumption \ref{a.3.2} is satisfied.
\begin{example} Consider equation (\ref{bsde}). Assume that:

\begin{longlist}[(a)]
\item[(a)]
$f(t,y,z)\dvtx[0,T]\times\mathbb{R}\times\mathbb{R}\rightarrow
\mathbb{R}$ is a deterministic function that has uniformly bounded
first- and second-order partial derivatives with respect to $y$ and
$z$, and $\int_0^Tf(t,0,0)^2\,dt<\infty$.

\item[(b)] The terminal value $\xi$ is a multiple stochastic
integral of the form
%
%e2.31 ###
%
\begin{equation}\label{mul}
\xi=\int_{[0,T]^n}g(t_1,\ldots,t_n)\,dW_{t_1}\cdots dW_{t_n},
\end{equation}
where
$n\geq2$ is an integer and $g(t_1,\ldots,t_n)$ is a symmetric
function in\break $L^2([0,T]^n)$, such that
\begin{eqnarray*}
\sup_{0\le u\le
T}\int_{[0,T]^{n-1}}g(t_1,\ldots,t_{n-1},u)^2\,dt_1\cdots
dt_{n-1}&<&\infty,
\\
\sup_{0\le u,v\le
T}\int_{[0,T]^{n-2}}g(t_1,\ldots,t_{n-2},u,v)^2\,dt_1\cdots
dt_{n-2}&<&\infty,
\end{eqnarray*}
and there exists a constant $L>0$ such that for any $u,v\in[0,T]$
\[
\int_{[0,T]^{n-1}}|g(t_1,\ldots,t_{n-1},u)-g(t_1,\ldots
,t_{n-1},v)|^2\,dt_1\cdots
dt_{n-1}<L|u-v|.
\]
\end{longlist}
From (\ref{mul}), we know that
\[
D_u\xi=n\int_{[0,T]^{n-1}}g(t_1,\ldots,t_{n-1},u)\,dW_{t_1}\cdots
dW_{t_{n-1}}.
\]
The above assumption implies Assumption \ref{a.3.2}, and
therefore, $Z$ satisfies the H\"{o}lder continuity property
(\ref{e-z}).
\end{example}
\begin{example}
Let $\Omega=C_0([0,1])$ equipped with the
Borel $\sigma$-field and Wiener measure.
Then, $\Omega$ is a Banach space with supremum norm \mbox{$\Vert\cdot\Vert
_\infty$},
and $W_t=\omega(t)$ is the canonical Wiener process. Consider equation
(\ref{bsde})
on the interval $[0,1]$. Assume that:
\begin{longlist}[(g3)]
\item[(g1)]
$f(t,y,z)\dvtx[0,1]\times\mathbb{R}\times\mathbb{R}\rightarrow
\mathbb{R}$
is a deterministic function that has uniformly\vspace*{1pt} bounded first- and
second-order partial derivatives with respect to $y$ and $z$, and
$\int_0^1f(t,0,0)^2\,dt<\infty$.\vspace*{1pt}
\item[(g2)] $\xi=\varphi(W)$, where
$\varphi\dvtx\Omega\to\mathbb{R}$ is twice Fr\'{e}chet
differentiable, and the first- and second-order Fr\'{e}chet
derivatives $\delta\varphi$ and $\delta^2\varphi$ satisfy
\[
\vert\varphi(\omega)\vert+\Vert\delta\varphi(\omega)\Vert
+\Vert\delta^2\varphi(\omega)\Vert\le
C_1\exp{\{C_2\Vert\omega\Vert_\infty^r\}}
\]
for all $\omega\in\Omega$ and some constants $C_1>0$, $C_2>0$ and
$0<r<2$, where \mbox{$\Vert\cdot\Vert$} denotes the operator norm (total
variation norm).

\item[(g3)] If $\lambda$ denotes the signed measure on $[0,1]$
associated with $\delta\varphi$, there exists a constant $L>0$
such that for all $0\le\theta\le\theta^\prime\le1$,
\[
\EE\vert\lambda((\theta,\theta^\prime])\vert^p\le
L\vert\theta-\theta^\prime\vert^{{p/2}}
\]
for some $p\geq2$.
\end{longlist}
It is easy\vspace*{1pt} to show that $D_\theta\xi=\lambda((\theta,1])$ and
$D_uD_\theta\xi=\nu((\theta,1]\times(u,1])$, where~$\nu$ denotes
the signed measure on $[0,1]\times[0,1]$ associated with
$\delta^2\varphi$. From the above assumptions and Fernique's
theorem, we can get Assumption~\ref{a.3.2}, and therefore, the
H\"{o}lder continuity property (\ref{e-z}) of $Z$.
\end{example}
\begin{example}\label{eg-2-11} Consider the following forward--backward
stochastic
differential equation (FBSDE for short):
%
%e2.32 ###
%
\begin{equation}\label{fbsde}
\cases{
\displaystyle X_t=X_0+\int_0^tb(r,X_r)\,dr+\int_0^t\sigma(r,X_r)\,dW_r,\vspace*{3pt}\cr
\displaystyle Y_t=\varphi\biggl(\int_0^T
X_r^2\,dr\biggr)+\int_t^Tf(r,X_r,Y_r,Z_r)\,dr-\int_t^TZ_r\,dW_r,}
\end{equation}
where $b, \sigma$, $\varphi$ and $f$ are deterministic functions,
and $X_0\in\mathbb{R}$.

We make the following assumptions:
\begin{longlist}[(h2)]
\item[(h1)] $b$ and $\sigma$ has uniformly bounded first- and
second-order partial derivatives with respect to $x$, and there is
a constant $L>0$, such that, for any $s, t\in[0,T]$,
$x\in\mathbb{R}$,
\[
\vert\sigma(t,x)-\sigma(s,x)\vert\le L\vert
t-s\vert^{1/2}.
\]
\item[(h2)] $\sup_{0\le t\le T}\{\vert
b(t,0)\vert+\vert\sigma(t,0)\vert\}<\infty$.
\item[(h3)] $\varphi$ is twice differentiable, and there exist a
constant $C>0$ and a~positive integer $n$ such that %
\[
\biggl\vert\varphi\biggl(\int_0^TX_t^2\,dt\biggr)\biggr\vert+
\biggl\vert\varphi'\biggl(\int_0^TX_t^2\,dt\biggr)\biggr\vert+\biggl\vert
\varphi''\biggl(\int_0^TX_t^2\,dt\biggr)\biggr\vert\le
C(1+\Vert X\Vert_\infty)^n,
\]
where $\Vert x\Vert_\infty=\sup\{|x(t)|, 0\le t\le T\}$ for any
$x\in C([0,T])$.
\item[(h4)] $f(t,x,y,z)$ has continuous and
uniformly bounded first- and second-order partial derivatives with
respect to $x, y$ and $z$ and \mbox{$\int_0^Tf(t,0,0,0)^2\,dt<\infty$}.
\end{longlist}
Notice that in this example,
$\Phi(X)=\varphi(\int_0^TX_t^2\,dt)$ is not necessarily
globally Lipschitz in $X$, and the results of \cite{Zh} cannot be
applied directly.

Under the above assumptions, (h1) and (h4), equation (\ref{fbsde}) has
a unique solution
triple $(X,Y,Z)$, and we have the following classical results: for any
real number $r>0$, there exists a constant $C>0$ such that%
\[
\EE\sup_{0\le t\le T}|X_t|^r<\infty,\qquad
\EE|X_t-X_s|^r\le C|t-s|^{{r/2}}
\]
for any $t,s\in[0,T]$.
For any fixed $(y,z)\in\mathbb{R}\times\mathbb{R}$, we have
$D_\theta f(t,X_t,y,z)=\partial_x f(t,X_t,y,z)D_\theta X_t$.
Then, under all the assumptions in this example, by Theorem 2.2.1 and
Lemma 2.2.2 in
\cite{N06} and the results listed above, we can verify Assumption \ref
{a.3.2}. Therefore,
$Z$ has the H\"{o}lder continuity property~(\ref{e-z}).

Note
that in the multidimensional case we do not require the matrix $\sigma
\sigma^T$ to be
invertible.
\end{example}
%

%s3 ###
\section{An explicit scheme for BSDEs}\label{sec3}

In the remaining part of this paper, we let $\pi=\{ 0=t_0<
t_1<\cdots< t_n=T\}$ be a partition of the interval $[0,T]$ and
$|\pi|=\max_{0\leq i\leq n-1}|t_{i+1}-t_i|$. Denote $\Delta
_i=t_{i+1}-t_i,%
0\leq i\leq n-1$.

From (\ref{bsde}), we know that, when $t\in
[t_i, t_{i+1}]$,
%
%e3.1 ###
%
\begin{equation} \label{e.4.1}
Y_t=Y_{t_{i+1}}+\int_t^{t_{i+1}}f(r,Y_r,Z_r)\,dr-\int_t^{t_{i+1}}Z_r\,dW_r.
\end{equation}
Comparing with the numerical schemes for forward stochastic
differential equations, we could introduce a numerical scheme of
the form
\begin{eqnarray}
Y^{1,\pi}_{t_n}&=&\xi^\pi,\nonumber\\
Y^{1,\pi}_{t_i}&=&Y_{t_{i+1}}^{1,\pi}
+f(t_{i+1},Y_{t_{i+1}}^{1,\pi} , Z^{1,\pi}_{t_{i+1}}
)\Delta_i -\int_{t_i}^{t_{i+1}}Z^{1,\pi}_r\,dW_r ,\nonumber\\
&&\eqntext{t\in[t_i,t_{i+1}) ,i=n-1,n-2,\ldots,0,}
\end{eqnarray}
where $\xi^\pi\in L^2(\Omega)$ is an approximation of the terminal
condition $\xi$. This leads to a backward recursive formula for the sequence
$\{Y^{1,\pi} _{t_i},Z^{1,\pi} _{t_i}\}_{0\le i\le n}$.
In fact, once $Y_{t_{i+1}}^{1,\pi}$ and $%
Z^{1,\pi}_{t_{i+1}}$ are defined, then we can find $Y_{t_{i}}^{1,\pi
}$ by
\[
Y_{t_{i}}^{1,\pi} =\mathbb{E}\bigl( Y_{t_{i+1}}^{1,\pi}
+f(t_{i+1},Y_{t_{i+1}}^{1,\pi} , Z^{1,\pi}_{t_{i+1}}
)\Delta_i %
|\mathcal{F}_{t_i}\bigr),
\]
and $\{Z^{1,\pi}_{r}\}_{t_i\le r<t_{i+1}}$ is determined by the
stochastic integral representation of the
random variable
\[
Y^{1,\pi}_{t_i}-Y_{t_{i+1}}^{1,\pi}
-f(t_{i+1},Y_{t_{i+1}}^{1,\pi} , Z^{1,\pi}_{t_{i+1}}
)\Delta_i.
\]
Although $\{Z^{1,\pi}_{r}\}_{t_i\le r<t_{i+1}}$ can be expressed
explicitly by
Clark--Ocone--Hauss\-man formula, its computation is a hard problem in practice.
On the other hand, there are difficulties in studying the convergence
of the above scheme.

An alternative scheme is
introduced in \cite{Zh}, where the approximating pairs
$(Y^\pi,Z^\pi)$ are defined recursively by
%
%e3.2 ###
%
\begin{eqnarray} \label{e.4.2}
Y^\pi_{t_n}&=&\xi^\pi,\qquad Z^\pi_{t_n}=0,\nonumber\\
Y^\pi_{t }&=&Y_{t_{i+1}}^\pi+f\biggl(t_{i+1},Y_{t_{i+1}}^\pi,\mathbb
{E}\biggl(\frac{1%
}{\Delta_{i+1}}\int_{t_{i+1}}^{t_{i+2}}Z^\pi_r\,dr \Big|\mathcal{F}%
_{t_{i+1}}\biggr)\biggr)\Delta_i \\
&&{} -\int_{t}^{t_{i+1}}Z^\pi_r\,dW_r ,\qquad
t\in[t_i,t_{i+1}) ,i=n-1,n-2,\ldots,0,\nonumber
\end{eqnarray}
where, by convention, $\mathbb{E}(\frac{1%
}{\Delta_{i+1}}\int_{t_{i+1}}^{t_{i+2}}Z^\pi_r\,dr |\mathcal{F}%
_{t_{i+1}})=0$ when $i=n-1$. In \cite{Zh} the following rate of
convergence is proved for this approximation scheme, assuming that the
terminal value $\xi$
and the generator $f$ are functionals of a~forward diffusion associated
with the BSDE,
%
%e3.3 ###
%
\begin{equation} \label{e.4.3}
\max_{0\leq i\leq n} \mathbb{E}| {Y}_{t_i}-Y_{t_i}^\pi|^2 +\mathbb
{E}%
\int_{0}^{T}| Z_t- Z_t^\pi|^2\,dt\leq{K}|\pi| .
\end{equation}

The main result of this section is the following, which on one hand
improves the above rate of convergence,
and on the other hand extends terminal
value $\xi$ and generator $f$ to more general situation.
\begin{theorem}
Consider the approximation scheme (\ref{e.4.2}). Let Assumption
\ref{a.3.2} be satisfied, and let the partition $\pi$ satisfy
\mbox{$\max_{0\le i\le n-1} {\Delta}%
_i/{\Delta}_{i+1}\le L_1$}, where $L_1$ is a constant. Assume that
a constant $L_2>0$ exists such that
%
%e3.4 ###
%
\begin{equation}\label{e-t}
|f(t_2, y, z)-f(t_1, y, z)|\le L_2 |t_2-t_1|^{1/2}
\end{equation}
for all $t_1, t_2\in[0,T]$ and $y, z\in\mathbb{R}$. Then there
are positive constants $K$ and~$\delta$, independent of the
partition $\pi$, such that, if $|\pi|<\delta$, then
%
%e3.5 ###
%
\begin{equation} \label{e.4.4}\quad
\mathbb{E}\sup_{0\le t\le T} | {Y}_{t}-Y_{t}^\pi|^2 +\mathbb{E}
\int_{0}^{T}| Z_t- Z_t^\pi|^2\,dt \leq
{K}(|\pi|+\EE|\xi-\xi^\pi|^2) .
\end{equation}
\end{theorem}
\begin{pf}
In this proof, $C>0$ will denote a constant independent of the
partition $\pi$, which may vary from line to line. Inequality
(\ref{e-z}) in Theorem~\ref{t.3.1}(b) yields the following
estimate (Theorem 3.1 in \cite{Zh}) with $p=2$:
\[
\sum_{i=0}^{n-1}\EE\int_{t_{i}}^{t_{i+1}}(\vert
Z_t-Z_{t_{i}}\vert^2+\vert Z_t-Z_{t_{i+1}}\vert^2)\,dt\le
C\vert\pi\vert.
\]
Using this estimate and following the same argument as the proof of Theorem
5.3 in \cite{Zh}, we can obtain the following result:
%
%e3.6 ###
%
\begin{equation}\label{e.4.5}
\max_{0\leq i\leq n}\mathbb{E} | {Y}_{t_i}-Y_{t_i}^\pi|^2
+\mathbb{E} \int_{0}^{T}| Z_t- Z_t^\pi|^2\,dt \leq
C(|\pi|+\EE|\xi-\xi^\pi|^2) .
\end{equation}
Denote
%
%e3.7 ###
%
\begin{equation}\label{e.4.10}
\widetilde Z_{t_i}^\pi=
\cases{
0, &\quad if $i=n$; \vspace*{2pt}\cr
\displaystyle \EE\biggl(\frac{1}{\Delta_{i }}\int_{t_{i
}}^{t_{i+1}}Z^\pi_r\,dr \Big|\mathcal{F}_{t_{i }}\biggr),
&\quad if $i=n-1,n-2,\ldots,0$.}
\end{equation}
If $t_i\le t<t_{i+1}$, $i=n-1,n-2,\ldots,0$, then, by iteration, we
have
%
%e3.8 ###
%
\begin{eqnarray}\label{e.4.11}
Y_t^\pi&=& Y_{t_{i+1}}^\pi+f(t_{i+1},Y_{t_{i+1}}^\pi
,\widetilde
Z_{t_{i+1}}^\pi)\Delta_i -\int_{t}^{t_{i+1}}Z^\pi_r\,dW_r
\nonumber\\[-8pt]\\[-8pt]
&=&\xi^\pi+\sum_{k=i+1}^{n } f(t_{k },Y_{t_{k }}^\pi,\tilde
Z_{t_{k }}^\pi)\Delta_{k-1} -\int_t^T Z^\pi_r\,dW_r .\nonumber
\end{eqnarray}
Therefore,
\[
Y_t^\pi=\mathbb{E}\Biggl( \xi^\pi+\sum_{k=i+1}^{n} f(t_{k
},Y_{t_{k }}^\pi,\widetilde Z_{t_{k
}}^\pi)\Delta_{k-1}\Big|\mathcal{F}_t\Biggr) ,\qquad t\in[t_i,t_{i+1}).
\]
We rewrite the BSDE
(\ref{bsde}) as follows:
%
%e3.9 ###
%
\begin{eqnarray} \label{e.4.12}
Y_{t} &=&\xi+\int_{t}^T f(r, Y_r, Z_r) \,dr-\int_{t}^T Z_r\,dW_r
\nonumber\\[-8pt]\\[-8pt]
&=&\xi+\sum_{k=i+1}^{n } f({t_{k }}, Y_{t_{k }}, Z_{t_{k }})
{\Delta}_{k-1} -\int_{t}^T Z_r\,dW_r+R_t^\pi,\nonumber
\end{eqnarray}
where
\begin{eqnarray*}
\vert R_t^\pi\vert&=&\Biggl\vert\int_{t}^{T} f(r, Y_r, Z_r)\,dr- \sum
_{k=i+1}^{n}f({%
t_{k }}, Y_{t_{k }}, Z_{t_{k }})\Delta_{k-1}\Biggr\vert\\
&=&\Biggl\vert\sum_{k=i+1}^n \int_{t_{k-1}}^{t_k} [ f(r, Y_r,
Z_r)-
f({ t_{k }}, Y_{t_{k }}, Z_{t_{k }})] \,dr-\int_{t_i}^tf(r, Y_r,
Z_r)\,dr\Biggr\vert\\
&\le& \sum_{k=i+1}^n \int_{t_{k-1}}^{t_k} | f(r, Y_r, Z_r)-
f({t_{k }}%
, Y_{t_{k }}, Z_{t_{k }})| \,dr +\int_{t_i}^{t_{i+1}}|f(r,
Y_r, Z_r)|\,dr.
\end{eqnarray*}
By Lemma \ref{T.2.1} and the Lipschitz condition on $f$, we have%
\[
\EE\biggl(\int_0^T|f(r,Y_r,Z_r)|^2\,dr\biggr)^{{p/2}}<\infty,
\]
and hence,%
%
%e3.10 ###
%
\begin{eqnarray}\label{e.4.12-1}
&&
\EE\max_{0\le i\le
n-1}\biggl(\int_{t_i}^{t_{i+1}}|f(r,Y_r,Z_r)|\,dr\biggr)^p\nonumber\\[-8pt]\\[-8pt]
&&\qquad\le|\pi
|^{{p/2}}\EE\biggl(\int_0^T|f(r,Y_r,Z_r)|^2\,dr\biggr)^{{p/2}}.\nonumber
\end{eqnarray}
Define a function $\{t(r)\}_{0\le r\le T}$ by
\[
t(r)=
\cases{
T, &\quad if $r=T$,
\cr
t_{i+1}, &\quad if $t_i\le r<t_{i+1}$, $i=n-1,\ldots,0$.}
\]
By the H\"{o}lder inequality, the boundedness of the first-order
partial derivatives of~$f$, (\ref{e-t}), (\ref{e-z}), Remark
\ref{r.3.8} and (\ref{e.4.12-1}), it is easy to see that
%
%e3.11 ###
%
\begin{eqnarray}\label{e.4.13}
\mathbb{E} \sup_{0\le t\le T} |R_t^\pi| ^p
&\le&2^{p-1}\biggl[ \mathbb{%
E}\biggl( \int_0^T \bigl| f(r, Y_r, Z_r)- f\bigl(t(r), Y_{t(r)},
Z_{t(r)}\bigr)\bigr|\,dr \biggr)^p \nonumber\\
&&\hspace*{67.5pt}{} +\EE\max_{0\le i\le
n-1}\biggl(\int_{t_i}^{t_{i+1}}|f(r,Y_r,Z_r)|\,dr\biggr)^p\biggr]
\nonumber\\
&\le& (2T)^{p-1} \mathbb{E}\int_0^T \bigl| f(r, Y_r, Z_r)-
f\bigl(t(r), Y_{t(r)}, Z_{t(r)}\bigr)\bigr|^p \,dr\\
&&{} +2^{p-1}|\pi|^{{p/2}}\EE\biggl(\int
_0^T|f(r,Y_r,Z_r)|^2\,dr\biggr)^{{p/2}} \nonumber\\
&\le& C |\pi|^{{p/2}},\nonumber
\end{eqnarray}
where, by convention, $R_T=0$.
In particular, we obtain
%
%e3.12 ###
%
\begin{equation} \label{e.4.14}
\mathbb{E} \sup_{0\le t\le T} | R_t^\pi|^2 \le C
|\pi| .
\end{equation}
To simplify the notation we denote
\[
\delta Y_t^\pi=Y_t-Y_t^\pi,\qquad \delta
Z_t^\pi=Z_t-Z_t^\pi\qquad \mbox{for all $t\in[0,T]$}
\]
and
\[
\widehat{Z}_{t_i}^\pi=Z_{t_i}-\widetilde Z_{t_i}^\pi
\qquad\mbox{for $i=n,n-1,\ldots,0$}.
\]
Then, when $t_i\le
t<t_{i+1}$, by (\ref{e.4.11}) and (\ref{e.4.12}) we can write
\begin{eqnarray*}
\delta Y_{t} ^\pi &=& \sum_{k=i+1}^n [ f({t_{k }}, Y_{t_{k }},
Z_{t_{k
}})- f({t_{k }}, Y_{t_{k }}^\pi, \widetilde Z_{t_{k }}^\pi)
]{\Delta}%
_{k-1}\\
&&{} -\int_{t}^T \delta Z_r^\pi\,dW_r+R_t^\pi+\delta\xi^\pi,
\end{eqnarray*}
where $\delta\xi^\pi=\xi-\xi^\pi$. Therefore, we obtain
%
%e3.13 ###
%
\begin{equation} \label{g1}
\delta Y_{t} ^\pi=\mathbb{E} \Biggl( \sum_{k=i+1}^n [
f({t_{k }},
Y_{t_{k }}, Z_{t_{k }})- f({t_{k }}, Y_{t_{k }}^\pi, \widetilde
Z_{t_{k }}^\pi)%
]{\Delta}_{k-1} +R_t^\pi+ \delta\xi^\pi\Big|\mathcal{F}%
_{t}\Biggr) .\hspace*{-36pt}
\end{equation}
Denote $\widetilde{f}_{t_k}^\pi=f(t_{k },Y_{t_{k }} , Z_{t_{k }})
-f(t_{k },Y_{t_{k }}^\pi,\widetilde Z_{t_{k }}^\pi)$. From equality
(\ref{g1}) for $t_j \le t < t_{j+1}$, where $i \le j \le n-1$, and
taking into account that
$\delta Y_T^\pi=\delta Y_{t_n}^\pi=\delta\xi^\pi$, we obtain
\[
% \label{e.4.16}
\sup_{t_i\le t\le T} |\delta Y_{t} ^\pi| \le\sup_{t_i\le t\le
T} \mathbb{E}\Biggl( \sum_{k=i+1}^{n }
|\tilde{f}_{t_k}^\pi|\Delta_{k-1}+ \sup_{0\le r\le T}
|R_r^\pi| +|\delta\xi^\pi|\Big|%
\mathcal{F}_{t}\Biggr) .
\]
The above conditional expectation is a martingale if it is
considered as a~process indexed by $t\in[t_i,T]$. Thus, using
Doob's maximal inequality, we obtain
\begin{eqnarray*} %\label{e.4.17}
\mathbb{E} \sup_{t_i\le t\le T} |\delta Y_{t}^\pi|^2
&\le&\mathbb{E} \sup_{t_i\le t\le T}\Biggl[ \mathbb{E}\Biggl(
\sum_{k=i+1}^{n } |\widetilde{f}_{t_k}^\pi
|\Delta_{k-1}+\sup_{0\le r\le T}| R_r^\pi|
+|\delta\xi^\pi| \Big|\mathcal{F}_{t}\Biggr)\Biggr]^2
\\
&\le& C \mathbb{E}\Biggl( \sum_{k=i+1}^{n }
|\widetilde{f}_{t_k}^\pi|\Delta_{k-1}+\sup_{0\le r\le T}|
R_r^\pi| +|\delta\xi^\pi| \Biggr)^2
\\
&\le& C \Biggl\{ \mathbb{E}\Biggl( \sum_{k=i+1}^{n }
|\widetilde{f}_{t_k}^\pi|\Delta_{k-1} \Biggr)^2 +
\mathbb{E}\sup_{0\le r\le T}| R_r^\pi|^2 +\EE|\delta
\xi^\pi|^2\Biggr\}.
\end{eqnarray*}
From (\ref{e.4.14}), we deduce
\[
\mathbb{E} \sup_{t_i\le t\le T} |\delta Y_{t}^\pi|^2 \le C
\Biggl\{ \mathbb{E}%
\Biggl( \sum_{k=i+1}^{n } |\widetilde{f}_{t_k}^\pi
|\Delta_{k-1} \Biggr)^2+\EE|\delta
\xi^\pi|^2 + |\pi|\Biggr\}.
\]
Using the Lipschitz condition on $f$, we obtain
%
%e3.14 ###
%
\begin{eqnarray} \label{e.4.19}
\mathbb{E} \sup_{t_i\le t\le T} |\delta Y_{t}^\pi|^2 &\le&
C\Biggl\{ (T-t_i)^2 \mathbb{E} \sup_{i+1\le k\le n} |\delta
Y_{t_k}^\pi|^2\nonumber\\[-2pt]
&&\hphantom{C\Biggl\{}
{} +
\mathbb{E} \Biggl(\sum_{k=i+1}^{n-1}
|\widehat{Z}_{t_k}^\pi|{\Delta}_{k-1}
\Biggr)^2 +\mathbb{E} |\widehat{Z}_{t_n}|^2{\Delta}_{n-1}^2\Biggr\}\\[-2pt]
&&{}
+C(\EE|\delta\xi^\pi|^2+|\pi|).\nonumber
\end{eqnarray}
Notice that
%
%e3.15 ###
%
\begin{eqnarray} \label{e.4.20}
\mathbb{E} \Biggl(\sum_{k=i+1}^{n-1}
|\widehat{Z}_{t_k}^\pi|{\Delta}_{k-1} \Biggr)^2
&=& \mathbb{E}
\Biggl( \sum_{k=i+1}^{n-1} \biggl| Z_{t_k}-\frac1{{\Delta}_k}
\int_{t_k}^{t_{k+1}} \mathbb{E} (Z_u^\pi|\mathcal{F}_{t_k})
\,du\biggr|%
{\Delta}_{k-1} \Biggr)^2 \nonumber\\[-2pt]
&\le& \mathbb{E} \Biggl( \sum_{k=i+1}^{n-1} \frac{{\Delta
}_{k-1}}{{\Delta}_k}
\int_{t_k}^{t_{k+1}} \mathbb{E} ( |Z_{t_k}-Z_u^\pi|
|\mathcal{F}%
_{t_k}) \,du\Biggr) ^2 \nonumber\\[-2pt]
&\le& L_1^2 \mathbb{E} \Biggl( \sum_{k=i+1}^{n-1}
\int_{t_k}^{t_{k+1}} \mathbb{E}
( |Z_{t_k}-Z_u^\pi||\mathcal{F}_{t_k}) \,du\Biggr)
^2 \nonumber\\[-10pt]\\[-10pt]
&\le& 2L_1^2\Biggl\{ \mathbb{E} \Biggl( \sum_{k=i+1}^{n-1}
\int_{t_k}^{t_{k+1}} \mathbb{E} (
|Z_{t_k}-Z_u||\mathcal{F}_{t_k}) \,du\Biggr) ^2
\nonumber\\[-2pt]
&&\hspace*{23pt}{} + \mathbb{E} \Biggl( \sum_{k=i+1}^{n-1} \int
_{t_k}^{t_{k+1}} \mathbb{%
E} ( |Z_u -Z_u^\pi||\mathcal{F}_{t_k}) \,du\Biggr)
^2\Biggr\}
\nonumber\\[-2pt]
&=&2L_1^2(I_1+I_2).\nonumber
\end{eqnarray}
Now the Minkowski and the H\"{o}lder inequalities yield
%
%e3.16 ###
%
\begin{eqnarray} \label{e.4.21}
I_1 &\le& \mathbb{E} \Biggl( \sum_{k=i+1}^{n-1} \biggl\{ \int
_{t_k}^{t_{k+1}}
\bigl( \mathbb{E} ( |Z_{t_k}-Z_u||\mathcal{F}_{t_k}
)%
\bigr)^2 \,du\biggr\}^{1/2} {\Delta}_k^{1/2} \Biggr) ^2 \nonumber\\[-2pt]
&\le& (T-t_{i }) \sum_{k=i+1}^{n-1} \int_{t_k}^{t_{k+1}}
\mathbb{E} \bigl( \mathbb{E} ( |Z_{t_k}-Z_u
||\mathcal{F}_{t_k})\bigr)^2 \,du
\nonumber\\[-10pt]\\[-10pt]
&\le& (T-t_{i }) \sum_{k=i+1}^{n-1} \int_{t_k}^{t_{k+1}}
\mathbb{E} |Z_{t_k}-Z_u|^2 \,du \nonumber\\[-2pt]
&\le&C(T-t_{i }) \sum_{k=i+1}^{n-1} \int_{t_k}^{t_{k+1}} |t_k-u|
\,du
\le C|\pi|.\nonumber
\end{eqnarray}
In a similar way and by (\ref{e.4.5}), we obtain
%
%e3.17 ###
%
\begin{eqnarray}\label{e.4.22}
I_2 &\le& (T-t_{i }) \sum_{k=i+1}^{n-1} \int_{t_k}^{t_{k+1}}
\mathbb{E} |Z_u -Z_{u}^\pi| ^2 \,du\nonumber\\[-8pt]\\[-8pt]
&=&(T-t_{i })\int_{t_{i+1}}^T\EE|\delta{Z}_u^\pi|^2\,du
\le C|\pi|.\nonumber
\end{eqnarray}
On the other hand,
%
%e3.18 ###
%
\begin{equation} \label{e.4.23}
\mathbb{E} (
\widehat{Z}_{t_n}^\pi{\Delta}_{n-1})^2=\mathbb{E} |Z_{t_n}|^2
|\Delta_{n-1}|^2\le C|\pi|^2 .
\end{equation}
From (\ref{e.4.19})--(\ref{e.4.23}), we have
%
%e3.19 ###
%
\begin{eqnarray} \label{e.4.24}
\mathbb{E} \sup_{t_i\le t\le T} |\delta Y_{t}^\pi|^2
&\le&
C_1 (T-t_{i})^2 \mathbb{E} \sup_{i+1\le k\le n} |\delta
Y_{t_k}^\pi|^2\nonumber\\[-8pt]\\[-8pt]
&&{}+C_2(\EE|\delta\xi^\pi|^2+|\pi|) ,\nonumber
\end{eqnarray}
where $C_1$ and $C_2$ are two positive constants independent of
the partition $\pi$.%

We can find a constant $\delta>0$ independent of the partition
$\pi$, such that $C_1(3\delta)^2<\frac{1}{2}$ and $T>2\delta$.
Denote $l=[\frac{T}{2\delta}]$ ($[x]$ means the greatest integer
no larger than $x$). Then $l\geq1$ is an integer independent of
the partition $\pi$. If $|\pi|<\delta$, then for the partition
$\pi$ we can choose $n-1>i_1>i_2>\cdots>i_l\geq0$, such that,
$T-2\delta\in(t_{i_1-1},t_{i_1}]$,
$T-4\delta\in(t_{i_2-1},t_{i_2}], \ldots, T-2\delta
l\in[0,t_{i_l}]$ (with $t_{-1}=0$).

For simplicity, we denote $t_{i_0}=T$ and $t_{i_{l+1}}=0$. Each
interval $[t_{i_{j+1}},t_{i_j}], j=0,1,\ldots,l$, has length less
than $3\delta$, that is, $|t_{i_j}-t_{i_{j+1}}|<3\delta$. On each
interval $[t_{i_{j+1}},t_{i_j}], j=0,1,\ldots,l$, we consider the
recursive formula (\ref{e.4.2}), and (\ref{e.4.24}) becomes
%
%e3.20 ###
%
\begin{eqnarray}\label{s1}
\mathbb{E} \sup_{t_{i_{j+1}}\le t\le t_{i_j}} |\delta Y_{t}^\pi|^2
&\le& C_1 (t_{i_j}-t_{i_{j+1}})^2 \mathbb{E} \sup_{i_{j+1}+1\le
k\le i_j} |\delta Y_{t_k}^\pi|^2\nonumber\\[-8pt]\\[-8pt]
&&{} +C_2(\EE|\delta
Y_{t_{i_j}}^\pi|^2+|\pi|) .\nonumber
\end{eqnarray}
Using (\ref{s1}), we can obtain inductively
%
%e3.21 ###
%
\begin{eqnarray}\label{e.4.25}
&&\mathbb{E} \sup_{t_{i_{j+1}}\le t\le t_{i_j}} |\delta Y_{t}^\pi|^2
\nonumber\\
&&\qquad\le C_1 (t_{i_j}-t_{i_{j+1}})^2 \mathbb{E} \sup_{i_{j+1}+1\le
k\le i_j} |\delta
Y_{t_k}^\pi|^2 +C_2(\EE|\delta Y_{t_{i_j}}^\pi|^2+|\pi|
)\nonumber\\
&&\qquad\le C_1(t_{i_j}-t_{i_{j+1}})^2\cdots C_1(t_{i_j}-t_{i_j-1})^2\EE
|\delta Y_{t_{i_j}}^\pi|^2\nonumber\\
&&\qquad\quad{} + C_2(\EE|\delta Y_{t_{i_j}}^\pi
|^2+|\pi|)\nonumber\\
&&\qquad\quad\hspace*{10pt}{}\times
\bigl(1+C_1(t_{i_j}-t_{i_{j+1}})^2+C_1(t_{i_j}-t_{i_{j+1}})^2C_1(t_{i_j}-t_{i_{j+1}+1})^2\nonumber\\
&&\qquad\quad\hspace*{26pt}{}+\cdots+C_1(t_{i_j}-t_{i_{j+1}})^2C_1(t_{i_j}-t_{i_{j+1}+1})^2\cdots
C_1(t_{i_j}-t_{i_j-1})^2\bigr)\\
&&\qquad\le(C_1(3\delta)^2)^{i_j-i_{j+1}}\EE|\delta
Y_{t_{i_j}}^\pi|^2\nonumber\\
&&\qquad\quad{}+C_2(\EE|\delta Y_{t_{i_j}}^\pi|^2+|\pi|)\nonumber\\
&&\qquad\quad\hspace*{10pt}{}\times
\bigl(1+C_1(3\delta)^2+(C_1(3\delta)^2)^2+\cdots+(C_1(3\delta
)^2)^{i_j-i_{j+1}}\bigr)\nonumber\\
&&\qquad\le\EE|\delta
Y_{t_{i_j}}^\pi|^2+\frac{C_2}{1-C_1(3\delta)^2}(\EE|\delta
Y_{t_{i_j}}^\pi|^2+|\pi|)\nonumber\\
&&\qquad\le\EE|\delta Y_{t_{i_j}}^\pi|^2+2C_2(\EE|\delta
Y_{t_{i_j}}^\pi|^2+|\pi|)\nonumber\\
&&\qquad=(2C_2+1)\EE|\delta Y_{t_{i_j}}^\pi|^2+2C_2|\pi|.\nonumber
\end{eqnarray}
By recurrence, we obtain
%
%e3.22 ###
%
\begin{eqnarray}\label{e.4.25-1}\qquad
&&\mathbb{E} \sup_{t_{i_{j+1}}\le t\le t_{i_j}} |\delta
Y_{t}^\pi|^2\nonumber\\
&&\qquad\le(2C_2+1)^{j+1}\EE|\delta\xi^\pi|^2+C_2|\pi|
\bigl(1+(2C_2+1)+\cdots+(2C_2+1)^j\bigr)\nonumber\\[-8pt]\\[-8pt]
&&\qquad\le(2C_2+1)^{l+1}\EE|\delta\xi^\pi|^2+C_2|\pi|
\bigl(1+(2C_2+1)+\cdots+(2C_2+1)^l\bigr)\nonumber\\
&&\qquad\le\frac{3(2C_2+1)^{l+1}}{2}(\EE|\delta\xi^\pi|^2+|\pi|).\nonumber
\end{eqnarray}
Therefore, taking $C=\frac{3(2C_2+1)^{l+1}}{2}$, we obtain
\[
\mathbb{E} \sup_{0\le t\le T} |\delta Y_{t}^\pi|^2\le\max_{0\le
j\le l}\mathbb{E} \sup_{t_{i_{j+1}}\le t\le t_{i_j}} |\delta
Y_{t}^\pi|^2\le C(|\pi|+\EE|\xi-\xi^\pi|^2).
\]
Combining the above estimate with (\ref{e.4.5}), we know that
there exists a constant $K>0$ independent of the partition $\pi$,
such that
\[
\mathbb{E} \sup_{0\le t\le T} |Y_t- Y_{t}^\pi|^2+\mathbb{E}
\int_{0}^{T}| Z_t- Z_t^\pi|^2\,dt \le
K(|\pi|+\EE|\xi-\xi^\pi|^2).
\]
\upqed\end{pf}
\begin{remark}
The numerical scheme introduced before, as other similar schemes,
involves the computation of conditional expectations with respect
to the $\sigma$-field $\mathcal{F}_{t_{i+1}}$. To implement this
scheme in practice we need to
approximate these conditional expectations. Some
work has been done to solve this problem, and we refer the reader
to the references \cite{BD,BT} and \cite{GLW}.
\end{remark}
%

%s4 ###
\section{An implicit scheme for BSDEs}\label{sec4}

In this section, we propose an implicit numerical scheme for the
BSDE (\ref{bsde}). Define the approximating\vadjust{\goodbreak} pairs $(Y^\pi,
Z^\pi)$ recursively by
%
%e4.1 ###
%
\begin{eqnarray} \label{e.5.1}
\quad Y^\pi_{t_n}&=&\xi^\pi, \nonumber\\
\quad Y^\pi_{t}&=&Y_{t_{i+1}}^\pi+f\biggl(t_{i+1},Y_{t_{i+1}}^\pi,\frac
{1}{%
\Delta_{i}}\int_{t_i}^{t_{i+1}}Z^\pi_r\,dr \biggr)\Delta_i
-\int_{t}^{t_{i+1}}Z^\pi_r\,dW_r , \\
\eqntext{t\in[t_i,t_{i+1}), i=n-1,n-2,\ldots,0,}
\end{eqnarray}
where the partition $\pi$ and $\Delta_i$, $i=n-1,\ldots,0$, are
defined in Section \ref{sec3}, and $\xi^\pi$ is an approximation of the
terminal value $\xi$. In this recursive formula (\ref{e.5.1}),
on each subinterval $[t_i, t_{i+1}), i=n-1,\ldots,0$, the
nonlinear ``generator'' $f$ contains the information of $Z^\pi$ on
the same interval. In this sense, this formula is different from
formula (\ref{e.4.2}),
and (\ref{e.5.1}) is an equation for
$\{(Y_t^\pi,Z_t^\pi)\}_{t_i\le t<t_{i+1}}$. When $|\pi|$ is
sufficiently small, the existence and uniqueness of the solution
to the above equation can be established. In fact, equation
(\ref{e.5.1}) is of the following form:
%
%e4.2 ###
%
\begin{equation} \label{e.5.2}
Y_t=\xi+ g\biggl(\int_a^b Z_r\,dr\biggr)-\int_t^b Z_r\,dW_r ,\qquad
t\in[a,b] \mbox{ and } 0\le a<b\le T.\hspace*{-28pt}
\end{equation}
For the BSDE (\ref{e.5.2}), we have the following theorem.
\begin{theorem}\label{t.5.1}
Let $0\le a<b\le T$ and $p\geq2$. Let $\xi$ be
$\mathcal{F}_b$-measurable and $\xi\in L^p(\Omega)$. If there exists
a constant $L>0$ such that
$g\dvtx(\Omega\times\mathbb{R},\mathcal{F}_b\otimes
\mathcal{B})\rightarrow(\mathbb{R},\mathcal{B})$ satisfies
\[
|g(z_1)-g(z_2)|\le L|z_1-z_2|
\]
for all $z_1,z_2\in\mathbb{R}$ and $g(0)\in L^p(\Omega)$, then
there is a constant $\delta(p,L)>0$, such that, when $b-a<
\delta(p,L)$, equation (\ref{e.5.2}) has a unique solution
$(Y,Z)\in S^p_{\mathcal{F}}([a,b])\times
H^p_{\mathcal{F}}([a,b])$.
\end{theorem}
\begin{pf}
We shall use the fixed point theorem for the mapping from
$H^p_{\mathcal{F}}([a,b])$ into $H^p_{\mathcal{F}}([a,b])$ which
maps $z$ to $Z$, where $(Y,Z)$ is the solution of the following
BSDE:
%
%e4.3 ###
%
\begin{equation}\label{e-ab}
Y_t=\xi+ g\biggl(\int_a^b z_r\,dr\biggr)-\int_t^b
Z_r\,dW_r,\qquad t\in[a,b].
\end{equation}
In fact, by the martingale representation theorem, there exist a
progressively measurable
process $Z= \{Z_t\}_{a\le t\le b}$ such that $\mathbb{E} \int_a^b
Z_t^2\,dt <\infty$
and
\[
\xi+ g\biggl(\int_a^b z_r\,dr\biggr)=\EE\biggl(\xi+ g\biggl(\int_a^b
z_r\,dr\biggr) \Big| \mathcal{F}_a \biggr)+\int_a^b Z_t\,dW_t.
\]
By the integrability properties of $\xi, g(0)$ and $z$, one can
show that $Z\in H^p_{\mathcal{F}}([a,b])$. Define
$Y_t=\EE(\xi+ g(\int_a^b
z_r\,dr)|\mathcal{F}_t), t\in[a,b]$. Then $(Y,Z)$
satisfies equation (\ref{e-ab}). Notice that $Y$ is a
martingale.
Then by the Lipschitz condition on $g$, the integrability of
$\xi, g(0)$ and $z$, and Doob's maximal inequality, we can
prove that $Y\in S^p_{\mathcal{F}}([a,b])$.\vadjust{\goodbreak}

Let $z^1, z^2$ be two elements in the Banach space
$H^p_{\mathcal{F}}([a,b])$, and let $(Y^1,Z^1)$, $(Y^2,Z^2)$ be the
associated solutions, that is,
\[
Y^i_t=\xi+ g\biggl(\int_a^b z^i_r\,dr\biggr)-\int_t^b
Z^i_r\,dW_r,\qquad t\in[a,b], i=1,2.
\]
Denote
\[
\bar Y=Y^1-Y^2 ,\qquad \bar Z=Z^1-Z^2 ,\qquad \bar z=z^1-z^2 .
\]
Then
%
%e4.4 ###
%
\begin{equation}\label{e.5.3}
\bar{Y}_t=g\biggl(\int_a^b z_r^1\,dr\biggr)-g\biggl(\int_a^b
z_r^2\,dr\biggr)-\int_t^b \bar{Z}_r \,dW_r
\end{equation}
for all $t\in[a,b]$. So
\[
\bar{Y}_t=\EE\biggl(g\biggl(\int_a^b z_r^1\,dr\biggr)-g\biggl(\int_a^b
z_r^2\,dr\biggr)\Big|\mathcal{F}_t\biggr)
\]
for all $t\in[a,b]$. Thus by Doob's maximal inequality, we have
%
%e4.5 ###
%
\begin{eqnarray}\label{e.5.4}
\mathbb{E}\sup_{a\le t\le b} |\bar{Y}_t|^p&=&\EE\sup_{a\le t\le b}
\biggl|\EE\biggl(g\biggl(\int_a^b z_r^1\,dr\biggr)-g\biggl(\int_a^b
z_r^2\,dr\biggr)\Big|\mathcal{F}_t\biggr)\biggr|^p\nonumber\\
&\le& C\mathbb{E} \biggl| g\biggl(\int_a^b
z_r^1\,dr\biggr)-g\biggl(\int_a^b z_r^2\,dr\biggr)\biggr|^p
\nonumber\\[-8pt]\\[-8pt]
&\le& C \mathbb{E} \biggl| \int_a^b z_r^1\,dr- \int_a^b z_r^2\,dr
\biggr|^p
\nonumber\\
&\le& C (b-a)^{{p/2}} \mathbb{E} \biggl(\int_a^b |\bar{z}_r|
^2 \,dr\biggr)^{{p/2}} ,\nonumber
\end{eqnarray}
where $C>0$ is a generic constant depending on $L$ and $p$, which
may vary from line to line. From (\ref{e.5.3}), it is
easy to see
\[
\bar{Y}_t=\bar{Y}_a+\int_a^t \bar{Z}_r \,dW_r
\]
for all $t\in[a,b]$. Therefore, by the Burkholder--Davis--Gundy
inequality\break and~(\ref{e.5.4}), we have
%
%e4.6 ###
%
\begin{eqnarray}\label{e.5.4-2}
\mathbb{E}\biggl( \int_a^b
|\bar{Z}_r|^2\,dr\biggr)^{{p/2}}&\le& C\mathbb{E} \sup_{a\le
t\le
b}\biggl| \int_a^t \bar{Z}_r \,dW_r\biggr|^p \nonumber\\
&\le& C\Bigl[ \mathbb{E}|\bar{Y}_a|^p+\mathbb{E}
\sup_{a\le
t\le b}|\bar{Y}_t|^p\Bigr] \\
&\le&C (b-a)^{{p/2}} \mathbb{E} \biggl(\int_a^b |\bar{z}_r|
^2 \,dr\biggr)^{{p/2}} ,\nonumber
\end{eqnarray}
that is,
\[
\Vert\bar{Z}\Vert_{H^p}\leq
C_1(b-a)^{1/2}\Vert\bar{z}\Vert_{H^p},
\]
where $C_1$ is a positive constant depending only on $L$ and $p$.

Take $\delta(p,L)=1/C_1^2$. It is obvious that the mapping is a
contraction when $b-a<\delta(p,L)$, and hence there exists a
unique solution $(Y,Z)\in S^p_{\mathcal{F}}([a,b])\times
H^p_{\mathcal{F}}([a,b])$ to the BSDE (\ref{e.5.2}).
\end{pf}

Now we begin to study the convergence of the scheme
(\ref{e.5.1}).
\begin{theorem}
\label{t.5.2} Let Assumption \ref{a.3.2} be satisfied, and let
$\pi$ be any partition. Assume that $\xi^\pi\in L^p(\Omega)$ and
there exists a constant $L_1>0$ such that, for all
$t_1, t_2\in[0,T]$,
\[
|f(t_2, y, z)-f(t_1, y, z)|\le L_1 |t_2-t_1|^{1/2}.
\]
Then, there are two positive constants $\delta$ and $K$
independent of the partition~$\pi$, such that, when
$|\pi|<\delta$, we have
\[
\mathbb{E} \sup_{0\leq t\leq T} | {Y}_{t}-Y_{t}^\pi|^p +\mathbb{E}
\biggl(\int_{0}^{T}| Z_t- Z_t^\pi|^2\,dt\biggr)^{{p/2}} \leq
K(|\pi|^{{p/2}}+\EE|\xi-\xi^\pi|^p) .
\]
\end{theorem}
\begin{pf}
If $|\pi|<\delta(p,L)$, where $\delta(p,L)$ is the constant in
Theorem \ref{t.5.1}, then Theorem \ref{t.5.1} guarantees the
existence and uniqueness of $(Y^\pi,Z^\pi)$. Denote, for
$i=n-1,n-2,\ldots,0$,
\[
\widetilde Z_{t_{i+1}}^\pi=\frac1{t_{i+1}-t_{t_i}}
\int_{t_{i}}^{t_{i+1}} Z_r^\pi\,dr.
\]
Notice that $\{\widetilde{Z}^\pi_{t_i},\}_{i=n-1,n-2,\ldots,0}$ here is
different from that in Section \ref{sec3}. Then
\begin{eqnarray*}
Y_{t_i}^\pi&=&Y_{t_{i+1}}^\pi+ f(t_{i+1}, Y_{t_{i+1}}^\pi,
\widetilde Z_{t_{i+1}}^\pi){\Delta}_i\\
&&{}-\int_{t_i}^{t_{i+1}}
Z_r^\pi\,dW_r,\qquad i=n-1,n-2,\ldots,0.
\end{eqnarray*}
Recursively, we obtain
\begin{eqnarray*}
Y_{t_i}^\pi&=&\xi^\pi+ \sum_{k=i+1}^n f(t_k, Y_{t_k}^\pi,
\widetilde Z_{t_k}^\pi){\Delta}_{k-1}\\
&&{}-\int_{t_i}^{T} Z_r^\pi
\,dW_r,\qquad i=n-1,n-2,\ldots,0.
\end{eqnarray*}
Denote
\[
\delta\xi^\pi=\xi-\xi^\pi,\qquad \delta Y_t^\pi=Y_t-Y_t^\pi
,\qquad
\delta Z_t^\pi=Z_t-Z_t^\pi,\qquad t\in[0,T],
\]
and
\[
\widehat{Z}_{t_i}^\pi=Z_{t_i}-\widetilde Z_{t_i}^\pi,\qquad
i=n-1,\ldots,0.\vadjust{\goodbreak}
\]

If $t\in[t_i ,t_{i+1})$, $i=n-1, n-2,\ldots,0$, then by
iteration, we have
%
%e4.7 ###
%
\begin{eqnarray}\label{e.5.5}
\delta Y_{t}^\pi&=&\delta\xi^\pi+ \sum_{k=i+1}^n
[f(t_k, Y_{t_k}, Z_{t_k})- f(t_k,
Y_{t_k}^\pi, \widetilde
Z_{t_k}^\pi)]{\Delta}_{k-1} \nonumber\\[-8pt]\\[-8pt]
&&{} -\int_{t_i}^{T} \delta{Z}_r^\pi\,dW_r+ R_{t}^\pi,\nonumber
\end{eqnarray}
where $R_t^\pi$ is exactly the same as that in Section \ref{sec3}.%

Denote $\widetilde{f}_{t_k}^\pi=f(t_{k },Y_{t_{k }} , Z_{t_{k }})
-f(t_{k },Y_{t_{k }}^\pi,\widetilde Z_{t_{k }}^\pi)$. Then for
$t\in[t_i,t_{i+1}), i\,{=}\,n-1,\break n-2,\ldots,0$, we have
%
%e4.8 ###
%
\begin{equation} \label{g3}
\delta Y_{t}^\pi=\mathbb{E} \Biggl(\delta\xi^\pi+ \sum_{k=i+1}^n
\widetilde{f}_{t_k}^\pi{\Delta}_{k-1} + R_{t}^\pi\Big|\mathcal{F}
_{t}\Biggr) .
\end{equation}
From equality
(\ref{g3}) for $t_j \le t < t_{j+1}$, where $i \le j \le n-1$, and
taking into account that
$\delta Y_T^\pi=\delta Y_{t_n}^\pi=\delta\xi^\pi$, we obtain
\[
% \label{e.4.16}
\sup_{t_i\le t\le T} |\delta Y_{t} ^\pi| \le\sup_{t_i\le t\le T}
\mathbb{E}\Biggl( \sum_{k=i+1}^{n }
|\widetilde{f}_{t_k}^\pi|\Delta_{k-1}+ \sup_{0\le r\le
T} |R_r^\pi| +|\delta\xi^\pi|\Big|%
\mathcal{F}_{t}\Biggr) .
\]
The above conditional expectation is a martingale if it is
considered as a~process indexed by $t$ for $t\in[t_i,T]$. Using
Doob's maximal inequality, (\ref{e.4.13}), and the Lipschitz
condition on $f$, we have
\begin{eqnarray*} %\label{e.4.17}
\hspace*{-2pt}&&\mathbb{E} \sup_{t_i\le t\le T} |\delta Y_{t}^\pi|^p\\
\hspace*{-2pt}&&\qquad\le\mathbb{E} \sup_{t_i\le t\le T}\Biggl[ \mathbb{E}\Biggl(
\sum_{k=i+1}^{n } |\widetilde{f}_{t_k}^\pi
|\Delta_{k-1}+\sup_{0\le r\le T}| R_r^\pi|
+|\delta\xi^\pi| \Big|\mathcal{F}_{t}\Biggr)\Biggr]^p
\\
\hspace*{-2pt}&&\qquad\le C \mathbb{E}\Biggl( \sum_{k=i+1}^{n }
|\widetilde{f}_{t_k}^\pi|\Delta_{k-1}+\sup_{0\le r\le T}|
R_r^\pi| +|\delta\xi^\pi| \Biggr)^p
\\
\hspace*{-2pt}&&\qquad\le C \Biggl\{ \mathbb{E}\Biggl( \sum_{k=i+1}^{n }
|\widetilde{f}_{t_k}^\pi|\Delta_{k-1} \Biggr)^p +
\mathbb{E}\sup_{0\le r\le T}| R_r^\pi|^p +\EE|\delta
\xi^\pi|^p\Biggr\} \\
\hspace*{-2pt}&&\qquad\le C\Biggl\{ \mathbb{E}\Biggl( \sum_{k=i+1}^{n } |\delta
Y_{t_k}^\pi|\Delta_{k-1} \Biggr)^p +\mathbb{E} \Biggl(
\sum_{k=i+1}^n |\widehat{Z}_{t_k}^\pi| {\Delta}_{k-1}
\Biggr)^p+ |\pi|^{{p/2}} +\EE|\delta\xi^\pi|^p\Biggr\}
\\
\hspace*{-2pt}&&\qquad\le C\Biggl\{ (T-t_i)^p\mathbb{E} \sup_{i+1\le k\le n}
|\delta Y_{t_k}^\pi|^p\\
\hspace*{-2pt}&&\qquad\quad\hphantom{C\Biggl\{}
{} +\mathbb{E} \Biggl(
\sum_{k=i+1}^n |\widehat{Z}_{t_k}^\pi| {\Delta}_{k-1}
\Biggr)^p+ |\pi|^{{p/2}} +\EE|\delta\xi^\pi|^p\Biggr\},
\end{eqnarray*}
where, and in the following, $C>0$ denotes a generic constant
independent of the partition $\pi$ and may vary from line to line.
On the other hand, we have, by the H\"{o}lder continuity of $Z$ given by
(\ref{e-z}),
\begin{eqnarray*} %\label{e.5.7}
&&
\mathbb{E} \Biggl( \sum_{k=i+1}^n |\widehat{Z}_{t_k}^\pi
| {%
\Delta}_{k-1} \Biggr)^p
\\
&&\qquad= \mathbb{E} \Biggl( \sum_{k=i+1}
^n \biggl|Z_{t_k}
-\frac1{{\Delta}_{k-1} } \int_{t_{k-1}}^{t_k} Z_r^\pi\,dr
\biggr|{\Delta}%
_{k-1}\Biggr)^p \\
&&\qquad\le \mathbb{E} \Biggl( \sum_{k=i+1} ^n \int_{t_{k-1}}^{t_k}
|Z_{t_k}-Z_r| \,dr +\sum_{k=i+1} ^n \int_{t_{k-1}}^{t_k}
|Z_r-Z_r^\pi| \,dr \Biggr)^p
\\
&&\qquad\le C|\pi|^{{p/2}} +2^{p-1}\mathbb{E} \biggl(
\int_{t_i}^{T} |Z_r-Z_r^\pi| \,dr
\biggr)^p \\
&&\qquad\le C|\pi|^{{p/2}}+ 2^{p-1} (T-t_i)^{{p/2}}
\mathbb{E} \biggl( \int_{t_i}^{T} |Z_r-Z_r^\pi
|^2 \,dr \biggr)^{{p/2}} \\
&&\qquad=C|\pi|^{{p/2}}+ 2^{p-1} (T-t_i)^{{p/2}} \mathbb{E}
\biggl( \int_{t_i}^{T} |\delta Z_r^\pi|^2 \,dr
\biggr)^{{p/2}}.
\end{eqnarray*}
Hence, we obtain
%
%e4.9 ###
%
\begin{eqnarray} \label{e.5.8}
&&\mathbb{E} \sup_{t_i\le t\le T} |\delta Y_{t}^\pi|^p
\nonumber\\
&&\qquad\le C_1\biggl\{
(T-t_i)^p \mathbb{E} \sup_{i+1 \le k\le n} |\delta Y_{t_k}
|^p\nonumber\\[-8pt]\\[-8pt]
&&\qquad\quad
\hphantom{C_1\biggl\{}
{} +
(T-t_i)^{{p/2}} \mathbb{E} \biggl( \int_{t_i}^{T} |\delta
Z_r^\pi|^2 \,dr
\biggr)^{{p/2}}\nonumber\\
&&\hspace*{119.2pt}{} +|\pi|^{{p/2}}+\EE|\delta\xi^\pi|^p\biggr\} ,\nonumber
\end{eqnarray}
where $C_1$ is a constant independent of the partition $\pi$.
By the Burkholder--Davis--Gundy inequality, we have
%
%e4.10 ###
%
\begin{equation}\label{e.5.8-1}
\mathbb{E} \biggl( \int_{t_i}^{T}|\delta Z_r^\pi|^2 \,dr
\biggr)^{{p/2}}\le c_p \mathbb{E} \biggl\vert\int_{t_i}^T
\delta Z_r^\pi\,dW_r\biggr\vert^p .
\end{equation}
From (\ref{e.5.5}), we obtain
%
%e4.11 ###
%
\begin{equation}\label{e.5.8-2}
\int_{t_i}^T \delta Z_r^\pi\,dW_r = \delta\xi^\pi+\sum_{k=i+1}^n
\widetilde{f}_{t_k}^\pi{\Delta}_{k-1}+R_{t_i}^\pi-\delta
Y_{t_i}^\pi.
\end{equation}
Thus, from (\ref{e.5.8-1}) and (\ref{e.5.8-2}), we obtain
\begin{eqnarray*}
&&\mathbb{E} \biggl( \int_{t_i}^{T}|\delta Z_r^\pi|^2 \,dr
\biggr)^{{p/2}}
\\[-8pt]
&&\qquad\le C_p\Biggl\{
\mathbb{E}\Biggl| \sum_{k=i+1}^n \widetilde{f}_{t_k}^\pi
{\Delta}_{k-1}\Biggr|^p+ \mathbb{E} |\delta\xi^\pi
|^p +
\mathbb{E} |R_{t_i}^\pi|^p+\mathbb{E}|\delta
Y_{t_i}^\pi|^p\Biggr\}.
\end{eqnarray*}
Similar to (\ref{e.5.8}), we have
\begin{eqnarray*}
&&\mathbb{E} \biggl( \int_{t_i}^{T}|\delta Z_r^\pi|^2 \,dr
\biggr)^{{p/2}}\\
&&\qquad\le C_2\biggl\{
(T-t_i)^p \mathbb{E} \sup_{i+1 \le k\le n} |\delta Y_{t_k} |^p\\
&&\hphantom{C_2\biggl\{}
\qquad\quad{} +
(T-t_i)^{{p/2}} \mathbb{E} \biggl( \int_{t_i}^{T} |\delta
Z_r^\pi|^2 \,dr
\biggr)^{{p/2}}+|\pi|^{{p/2}}+\EE|\delta\xi^\pi|^p \biggr\},
\end{eqnarray*}
where $C_2$ is a constant independent of the partition $\pi$.

If $C_2(T-t_i)^{{p/2}}<\frac{1}{2}$, then we have
%
%e4.12 ###
%
\begin{eqnarray}\label{e.5.9}
\mathbb{E} \biggl( \int_{t_i}^{T}|\delta Z_r^\pi|^2 \,dr
\biggr)^{{p/2}} &\le& 2C_2(T-t_i)^p\mathbb{E} \sup_{i+1 \le
k\le n} |\delta Y_{t_k}
|^p\nonumber\\[-8pt]\\[-8pt]
&&{}+2C_2(|\pi|^{{p/2}}+\EE|\delta\xi^\pi|^p).\nonumber
\end{eqnarray}
Substituting (\ref{e.5.9}) into (\ref{e.5.8}), we have
%
%e4.13 ###
%
\begin{eqnarray}\label{e.5.10}
&&\mathbb{E} \sup_{t_i\le t\le T} |\delta Y_{t}^\pi|^p\nonumber\\
&&\qquad\le C_1\bigl(1+2C_2(T-t_i)^{{p/2}}\bigr)
(T-t_i)^p \mathbb{E} \sup_{i+1 \le k\le n} |\delta Y_{t_k} |^p
\nonumber\\[-8pt]\\[-8pt]
&&\qquad\quad{} +
C_1\bigl(1+2C_2(T-t_i)^{{p/2}}\bigr)(|\pi|^{
{p}/{2}}+\EE|\delta\xi^\pi|^p)\nonumber\\
&&\qquad\le 2C_1 (T-t_i)^p \mathbb{E} \sup_{i+1 \le k\le n} |\delta
Y_{t_k}
|^p+2C_1(|\pi|^{{p/2}}+\EE|\delta\xi^\pi|^p) .\nonumber
\end{eqnarray}
We can find a positive constant $\delta<\delta(p,L)$ independent
of the partition $\pi$, such that,
%
%e4.15 ###
%e4.14 ###
%
\begin{eqnarray}
\label{s10}
C_2(3\delta)^{{p/2}}&<&\tfrac{1}{2},\\
\label{s11}
2C_1(3\delta)^p&<&\tfrac{1}{2}
\end{eqnarray}
and $T>2\delta$. Denote $l=[\frac{T}{2\delta}]$. Then $l\geq1$ is
an integer independent of the partition $\pi$. If $|\pi|<\delta$,
then for the partition $\pi$ we can choose
$n-1>i_1>i_2>\cdots>i_l\geq0$, such that,
$T-2\delta\in(t_{i_1-1},t_{i_1}]$,
$T-4\delta\in(t_{i_2-1},t_{i_2}], \ldots, T-2\delta
l\in[0,t_{i_l}]$ (with $t_{-1}=0$). For simplicity, we denote
$t_{i_0}=T$ and $t_{i_{l+1}}=0$. Each interval
$[t_{i_{j+1}},t_{i_j}], j=0,1,\ldots,l$, has length less than
$3\delta$, that is, $|t_{i_j}-t_{i_{j+1}}|<3\delta$. On
$[t_{i_{j+1}},t_{i_j}]$, we consider the recursive formula
(\ref{e.5.1}). Then (\ref{e.5.10})--(\ref{s11}) yield
%
%e4.16 ###
%
\begin{eqnarray}\label{e.5.11}\quad
&&\mathbb{E} \sup_{t_{i_{j+1}}\le t\le t_{i_{j}}} |\delta Y_{t}^\pi
|^p\nonumber\\
&&\qquad\le 2C_1(t_{i_{j}}-t_{i_{j+1}})^p \mathbb{E} \sup_{i_{j+1}+1 \le
k\le i_j}
|\delta Y_{t_k} |^p+ 2C_1(|\pi|^{{p/2}}+\EE|\delta
Y_{t_{i_j}}^\pi|^p)\nonumber\\[-8pt]\\[-8pt]
&&\qquad\le 2C_1(3\delta)^p \mathbb{E} \sup_{i_{j+1}+1 \le k\le i_j}
|\delta Y_{t_k} |^p+
2C_1(|\pi|^{{p/2}}+\EE|\delta Y_{t_{i_j}}^\pi|^p
)\nonumber\\
&&\qquad\le\frac{1}{2}\sup_{i_{j+1}+1 \le k\le i_j} |\delta Y_{t_k}
|^p+2C_1(|\pi|^{{p/2}}+\EE|\delta
Y_{t_{i_j}}^\pi|^p).\nonumber
\end{eqnarray}
As in the proof of (\ref{e.4.25}) and (\ref{e.4.25-1}), we have
\[
\mathbb{E} \sup_{t_{i_{j+1}}\le t\le t_{i_{j}}} |\delta Y_{t}^\pi
|^p \le(4C_1+1)\EE|\delta Y_{t_{i_j}}^\pi
|^p+4C_1|\pi|^{{p/2}}
\]
and
\[
\mathbb{E} \sup_{t_{i_{j+1}}\le t\le t_{i_{j}}} |\delta Y_{t}^\pi
|^p
\le\frac{3(4C_1+1)^{l+1}}{2}(\EE|\delta\xi^\pi|^2+|\pi
|^{{p/2}}).
\]
Therefore, we obtain
%
%e4.17 ###
%
\begin{eqnarray}\label{e.5.12}
\mathbb{E} \sup_{0\le t\le T} |\delta Y_{t}^\pi|^p&\le&\max_{0\le
j\le l}\mathbb{E} \sup_{t_{i_{j+1}}\le t\le t_{i_{j}}} |\delta
Y_{t}^\pi|^p\nonumber\\[-8pt]\\[-8pt]
&\le&\frac{3(4C_1+1)^{l+1}}{2}(\EE|\delta\xi^\pi
|^p+|\pi|^{{p/2}}).\nonumber
\end{eqnarray}
On $[t_{i_{j+1}},t_{i_j}], j=0,1,\ldots,l$, based on the recursive
formula (\ref{e.5.1}) and (\ref{e.5.12}), inequality (\ref{e.5.9})
becomes
\begin{eqnarray*}
&&\mathbb{E} \biggl( \int_{t_{i_{j+1}}}^{t_{i_{j}}}|\delta Z_r^\pi
|^2 \,dr
\biggr)^{{p/2}}
\\
&&\qquad\le 2C_2(t_{i_{j}}-t_{i_{j+1}})^p\mathbb{E} \sup_{i_{j+1}+1 \le
k\le i_j} |\delta Y_{t_k}
|^p+2C_2(|\pi|^{{p/2}}+\EE|\delta\xi^\pi|^p)\\
&&\qquad\le 2C_2(3\delta)^p\mathbb{E} \sup_{i_{j+1}+1 \le k\le i_j}
|\delta Y_{t_k}
|^p+2C_2(|\pi|^{{p/2}}+\EE|\delta\xi^\pi|^p)\\
&&\qquad\le \frac{1}{2}\mathbb{E} \sup_{i_{j+1}+1 \le k\le i_j} |\delta
Y_{t_k}
|^p+2C_2(|\pi|^{{p/2}}+\EE|\delta\xi^\pi|^p)
\\
&&\qquad\le \biggl(\frac{3(4C_1+1)^{l+1}}{4}+2C_2\biggr)(|\pi
|^{{p/2}}+\EE|\delta\xi^\pi|^p).
\end{eqnarray*}
Thus
%
%e4.18 ###
%
\begin{eqnarray}\label{e.5.13}\qquad
&&\EE\biggl(\int_0^T|\delta Z_t^\pi|^2\,dt\biggr)^{
{p/2}}\nonumber\\
&&\qquad=\EE\Biggl(\sum_{j=0}^{l}\int_{t_{i_{j+1}}}^{t_{i_{j}}}|\delta
Z_t^\pi
|^2 \,dt\Biggr)^{{p/2}}\nonumber\\[-8pt]\\[-8pt]
&&\qquad\le(l+1)^{{p}/{2}-1}\sum_{j=0}^{l}\mathbb{E} \biggl(
\int_{t_{i_{j+1}}}^{t_{i_{j}}}|\delta Z_t^\pi|^2 \,dt
\biggr)^{{p/2}}\nonumber\\
&&\qquad\le(l+1)^{{p/2}}\biggl(\frac{3(4C_1+1)^{l+1}}{4}+2C_2
\biggr)(|\pi|^{{p/2}}+\EE|\delta\xi^\pi|^p).\nonumber
\end{eqnarray}
Combining (\ref{e.5.12}) and (\ref{e.5.13}), we know that there
exists a constant
\[
K=(l+1)^{{p/2}}\biggl(\frac{3(4C_1+1)^{l+1}}{2}+4C_2\biggr)
\]
independent of the partition $\pi$, such that
\begin{eqnarray*}
&&\mathbb{E} \sup_{0\leq t\leq T} | {Y}_{t}-Y_{t}^\pi|^p +\mathbb{E}
\biggl(\int_{0}^{T}| Z_t- Z_t^\pi|^2\,dt\biggr)^{{p/2}} \\
&&\qquad\leq
K(|\pi|^{{p/2}}+\EE|\xi-\xi^\pi|^p) .
\end{eqnarray*}
\upqed\end{pf}
\begin{remark}
The advantages of this implicit numerical scheme are:

\begin{longlist}
\item we can obtain the rate of convergence in $L^p$
sense;%
\item the partition $\pi$ can be arbitrary ($|\pi|$
should be small enough) without assuming $\max_{0\le i\le n-1} {\Delta
}%
_i/{\Delta}_{i+1}\le L_1$.
\end{longlist}
\end{remark}

%s5 ###
\section{A new discrete scheme}\label{sec5}

For all the numerical schemes considered in Sections \ref{sec3} and \ref
{sec4}, one
needs to evaluate processes $\{Z_t^\pi\}_{0\le t\le T}$ with
continuous index $t$. In this section, we use the representation
of $Z$ in terms of the Malliavin derivative of $Y$ to derive a
completely discrete scheme.

From (\ref{e.3.12}), $\{D_\theta Y_t\}_{0\le\theta\le
t\le T}$ can be represented as
%
%e5.1 ###
%
\begin{equation}\label{e.6.3}
D_{\theta} Y_t=\mathbb{E} \biggl( \rho_{t, T} D_{\theta} \xi+\int_t^T
\rho_{t, r } D_{\theta} f(r, Y_r, Z_r) \,dr\Big|\mathcal{F}_t\biggr),
\end{equation}
where
%
%e5.2 ###
%
\begin{equation}\label{e.6.4}
\rho_{t, r} =\exp\biggl\{ \int_t^r \beta_s\,dW_s +\int_t^r
\biggl( \alpha_s -\frac12 \beta_s^2%
\biggr)\,ds\biggr\}
\end{equation}
with $\alpha_s=\partial_y f(s, Y_s, Z_s)$ and $\beta_s=\partial_
z f(s, Y_s, Z_s)$.

Using that $Z_t=D_tY_t$, $\mu\times P$ a.e., from (\ref
{bsde}), %
(\ref{e.6.3}) and (\ref{e.6.4}), we propose the
following numerical scheme.
We define
recursively
%
%e5.3 ###
%
\begin{eqnarray}\label{e.6.5}\qquad
Y_{t_n}^\pi&=&\xi,\qquad Z_{t_n}^\pi=D_{T}\xi,\nonumber\\
Y_{t_i}^\pi
&=&\mathbb{E}\bigl( Y_{t_{i+1}}^\pi+f(t_{i+1}, Y_{t_{i+1}}^\pi,
Z_{t_{i+1}}^\pi){\Delta}_i |\mathcal{F}_{t_i}\bigr) ,
\nonumber\\[-9pt]\\[-9pt]
Z_{t_i} ^\pi&=&\mathbb{E} \Biggl( \rho_{t_{i+1} , t_n}^\pi D_{t_i}
\xi+\sum_{k=i}^{n-1} \rho_{t_{i+1}, t_{k+1}}^\pi D_{t_i}
f(t_{k+1}, Y_{t_{k+1}}^\pi, Z_{t_{k+1}}^\pi) {\Delta} _k
\Big|\mathcal{F}_{t_i}\Biggr),\nonumber\\
\eqntext{i=n-1, n-2,\ldots, 0,}
\end{eqnarray}
where $\rho_{t_i, t_i}^\pi=1, i=0,1,\ldots,n$, and for $0\le i<
j\le n$,
%
%e5.4 ###
%
\begin{eqnarray} \label{e.6.7}\qquad
\rho_{t_i, t_j}^\pi&=& \exp\Biggl\{ \sum_{k=i}^{j-1}
\int_{t_k}^{t_{k+1} }
\partial_ z f(r, Y_{t_k}^\pi, Z_{t_k}^\pi)\,dW_r \nonumber\\[-9pt]\\[-9pt]
&&\hphantom{\exp\Biggl\{}{} +\sum_{k=i}^{j-1} \int_{t_k}^{t_{k+1} } \biggl( \partial_y f(r,
Y_{t_k}^\pi, Z_{t_k}^\pi) -\frac12 [\partial_ z f(r, Y_{t_k}^\pi,
Z_{t_k}^\pi)]^2\biggr)\,dr \Biggr\} .\nonumber
\end{eqnarray}
An alternative expression for $\rho_{t_i, t_j}^\pi$ is given
by the following formula:
%
%e5.5 ###
%
\begin{eqnarray} \label{e.6.8}
\rho_{t_i, t_j}^\pi&=& \exp\Biggl\{ \sum_{k=i}^{j-1} \partial_ z
f(t_k,
Y_{t_k}^\pi, Z_{t_k}^\pi)(W_{t_{k+1}}-W_{t_k}) \nonumber\\[-9pt]\\[-9pt]
&&\hphantom{\exp\Biggl\{}
{} +\sum_{k=i}^{j-1} \biggl( \partial_y f({t_k}, Y_{t_k}^\pi,
Z_{t_k}^\pi) -
\frac12 [\partial_ z f({t_k}, Y_{t_k}^\pi, Z_{t_k}^\pi)]^2
\biggr){\Delta}_k %
\Biggr\}.\nonumber
\end{eqnarray}
However, we will only consider the scheme (\ref{e.6.5})
with $\rho_{t_i, t_j}^\pi$ given by (\ref{e.6.7}).

We make the following assumptions:

\begin{longlist}[(G2)]
\item[(G1)] $f(t,y,z)$ is deterministic, which implies
$D_{\theta} f(t,y,z)=0$.

\item[(G2)] $f(t,y,z)$ is linear with respect to $y$ and $z$;
namely, there are three functions $g(t)$, $h(t)$ and $f_1(t)$ such that
\[
f(t,y,z)=g(t)y+h(t) z+f_1(t) .
\]
Assume that $g$, $h$ are bounded and $f_1\in L^2([0,T])$.
Moreover, there exists a~constant $L_2>0$, such that, for all
$t_1, t_2\in[0,T]$,
\[
|g(t_2)-g(t_1)|+|h(t_2)-h(t_1)|+|f_1(t_2)-f_1(t_1)|\le L|t_2-t_1|^{1/2}.
\]

\item[(G3)] $ \EE\sup_{0\le\theta\le T}\vert
D_\theta\xi\vert^r<\infty, $ for all $r\geq1$.
\end{longlist}
Notice that (G1) and (G2) imply (ii) and (iii) in Assumption
\ref{a.3.2}.
\begin{remark}
We propose condition (G1) in order to simplify\break
$\{Z_{t_i}^\pi\}_{i=n-1,\ldots,0}$ in formula (\ref{e.6.5}). In
fact, there are some difficulties in generalizing the condition
(G)s, especially (G1), to a forward--backward stochastic
differential equation (\mbox{FBSDE}, for short) case.\vadjust{\goodbreak}

If we consider a FBSDE
\[
\cases{
\displaystyle X_t=X_0+\int_0^tb(r,X_r)\,dr+\int_0^t\sigma(r,X_r)\,dW_r,\vspace*{2pt}\cr
\displaystyle Y_t=\xi+\int_t^Tf(r,X_r,Y_r,Z_r)\,dr-\int_t^TZ_r\,dW_r,}
\]
where $X_0\in\mathbb{R}$, and the functions $b, \sigma, f$ are
deterministic, then under some appropriate conditions [e.g., (h1)--(h4)
in Example \ref{eg-2-11}] $Z^\pi_{t_i}$ for
$i=n-1,\ldots,0$ in (\ref{e.6.5}) is of the form
\begin{eqnarray*}
Z_{t_i} ^\pi&=&\mathbb{E} \Biggl( \rho_{t_{i+1} , t_n}^\pi D_{t_i}
\xi\\[-3pt]
&&\hphantom{\mathbb{E} \Biggl(}
{}+\sum_{k=i}^{n-1} \rho_{t_{i+1}, t_{k+1}}^\pi
\partial_x f(t_{k+1},X_{t_{k+1}}^\pi, Y_{t_{k+1}}^\pi,
Z_{t_{k+1}}^\pi)D_{t_i}X_{t_{k+1}}^\pi{\Delta} _k
\Big|\mathcal{F}_{t_i}\Biggr),
\end{eqnarray*}
where $(X^\pi, Y^\pi,Z^\pi)$ is a certain numerical scheme for
$(X,Y,Z)$. It is hard to guarantee the existence and the
convergence of Malliavin derivative of $X^\pi$, and therefore, the
convergence of $Z^\pi$ is difficult to derive.
\end{remark}
\begin{theorem}\label{t.6.1}
Let Assumption \ref{a.3.2}\textup{(i)} and assumptions \textup{(G1)--(G3)} be
satisfied. Then there are positive constants $K$ and $\delta$
independent of the partition $\pi$, such that, when $|\pi|<\delta$
we have
\[
\mathbb{E} \max_{0\le i\le n} \{ |Y_{t_i}-Y_{t_i}^\pi|^p +
|Z_{t_i}-Z_{t_i}^\pi|^p\}\le K |\pi|^{{p}/{2}-
{p}/({2\log({1}/{|\pi|})%
})} \biggl( \log\frac{1}{|\pi|}\biggr)^{{p/2}} .
\]
\end{theorem}
\begin{pf} In the proof, $C>0$ will denote a constant
independent of the partition~$\pi$, which may vary from line to
line. Under the assumption (G1), we can see that
%
%e5.6 ###
%
\begin{equation}
Z_{t_i}^\pi=\mathbb{E} ( \rho_{t_{i+1}, t_n}^\pi D_{t_i} \xi
|%
\mathcal{F}_{t_i}),\qquad i=n-1,n-2,\ldots,0.
\end{equation}
Denote, for $i=n-1,n-2,\ldots,0$,
\[
\delta Z_{t_i}^\pi=Z_{t_i} -Z_{t_i}^\pi, \qquad\delta
Y_{t_i}^\pi=Y_{t_i} -Y_{t_i}^\pi.\vadjust{\goodbreak}
\]
Since $|e^x-e^y|\le(e^x+e^y)|x-y|$, we deduce, for all
$i=n-1,n-2,\ldots,0$,
\begin{eqnarray*}
|\delta Z_{t_i}^\pi| &=& \bigl|\mathbb{E} ( \rho_{t_{i}, t_n}
D_{t_i} \xi%
|\mathcal{F}_{t_i})-\mathbb{E} ( \rho_{t_{i+1},
t_n}^\pi
D_{t_i} \xi|\mathcal{F}_{t_i})\bigr| \\
&\le&\mathbb{E} ( | \rho_{t_{i}, t_n} - \rho_{t_{i+1},
t_n}^\pi| |D_{t_i} \xi||\mathcal{F}_{t_i} ) \\
&\le& \mathbb{E} \biggl( |D_{t_i} \xi| (
\rho_{t_{i}, t_n} + \rho_{t_{i+1}, t_n}^\pi)\\
&&\hphantom{\mathbb{E} \biggl(}
{}\times \biggl|
\int_{t_i}^T h(r) \,dW_r + \int_{t_i}^T
g(r) \,dr
-\frac12 \int_{t_i}^T h(r)^2 \,dr\\
&&\hphantom{\mathbb{E} \biggl({}\times \biggl|}
{}-\sum_{k=i+1}^{n-1}
\int_{t_k}^{t_{k+1}} h(r) \,dW_r- \sum_{k=i+1}^{n-1}
\int_{t_k}^{t_{k+1}}
g(r) \,dr \\
&&\hspace*{123pt}{} +\frac12\sum_{k=i+1}^{n-1} \int_{t_k}^{t_{k+1}} h(r)^2
\,dr\biggr|\Big| %
\mathcal{F}_{t_i} \biggr) \\
&\le& \mathbb{E} \biggl( |D_{t_i} \xi| (
\rho_{t_{i}, t_n} + \rho_{t_{i+1}, t_n}^\pi) \\
&&\hphantom{\mathbb{E} \biggl(}
{}\times\biggl[
\biggl|\int_{t_i}^{t_{i+1}} h(r) \,dW_r\biggr| +\int_{t_i}^{t_{i+1}}
|g(r)|\,dr\\
&&\hphantom{\mathbb{E} \biggl({}\times\biggl[}\hspace*{45pt}
{} +\frac12 \int_{t_i}^{t_{i+1}} h(r)^2 \,dr
\biggr]\Big|\mathcal{F}_{t_i} \biggr) .
\end{eqnarray*}
From (G2), we have
\begin{eqnarray*}
&&|D_{t_i} \xi| \rho_{t_{i+1}, t_n}^\pi\\
&&\qquad\le |D_{t_i} \xi|
\exp\Biggl\{ \int_{t_{i+1}}^T h(r) \,dW_r +\sum_{k=i+1}^{n-1}
\int_{t_k}^{t_{k+1}} g(r)\,dr-\frac12 \int_{t_{i+1}}^T h(r)^2 \,dr
\Biggr\} \\
&&\qquad\le C_1\Bigl({\sup_{0\le\theta\le
T}}|D_\theta\xi|\Bigr)\biggl(\sup_{0\le t\le
T}\exp\biggl\{\int_t^Th(r)\,dW_r\biggr\}\biggr),
\end{eqnarray*}
where $C_1>0$ is a constant independent of the partition $\pi$.

In the same way, we obtain
\[
|D_{t_i} \xi|\rho_{{t_i},t_n}<C_1\Bigl({\sup_{0\le\theta\le
T}}|D_\theta\xi|\Bigr)\biggl(\sup_{0\le t\le
T}\exp\biggl\{\int_t^Th(r)\,dW_r\biggr\}\biggr) .
\]
Thus for $i=n-1,n-2,\ldots,0$,
\begin{eqnarray*}
|\delta Z_{t_i}^\pi|
&\le& 2C_1 \mathbb{E}
\biggl(\Bigl({\sup_{0\le\theta\le
T}}|D_\theta\xi|\Bigr)\biggl(\sup_{0\le t\le
T}\exp\biggl\{\int_t^Th(r)\,dW_r\biggr\}\biggr)\\
&&\hphantom{2C_1 \mathbb{E}
\biggl(}
{}\times\biggl[
\biggl|\int_{t_i}^{t_{i+1}} h(r) \,dW_r\biggr|
+\int_{t_i}^{t_{i+1}} |g(r)|\,dr +\frac12 \int_{t_i}^{t_{i+1}}
h(r)^2 \,dr
\biggr]\Big|\mathcal{F}_{t_i} \biggr)\\
&\le& 2C_1 \mathbb{E} \biggl( \Bigl({\sup_{0\le\theta\le
T}}|D_\theta\xi|\Bigr)\biggl(\sup_{0\le t\le
T}\exp\biggl\{\int_t^Th(r)\,dW_r\biggr\}\biggr)\\
&&\hphantom{2C_1 \mathbb{E} \biggl(}
{}\times
\biggl[{\sup_{0\le
k\le n-1}} \biggl|\int_{t_k}^{t_{k+1}} h(r)
\,dW_r\biggr|
+ \sup_{0\le k\le n-1}\int_{t_k}^{t_{k+1}} |g(r)|\,dr\\
&&\hspace*{153.8pt}{} +\frac12
\sup_{0\le k\le n-1} \int_{t_k}^{t_{k+1}} h(r)^2 \,dr\biggr]
\Big|\mathcal{F}_{t_i} %
\biggr) .
\end{eqnarray*}
The right-hand side of the above inequality is a martingale as a
process indexed by $i=n-1,n-2,\ldots,0$.

Let $\eta_t=\exp\{-\int_0^th(u)\,dW_u\}$. Then, $\eta_t$
satisfies the following linear stochastic differential equation:
\[
\cases{
d\eta_t=-h(t)\eta_t\,dW_t+\frac{1}{2}h(t)^2\eta_tdt,\cr
\eta_0=1.}
\]
By (G1), (G2), the H\"{o}lder inequality and Lemma \ref{l.3.1},
it is easy to show that, for any $r\geq0$,
%
%e5.7 ###
%
\begin{eqnarray}\label{ee.6.1}\quad
&&\EE\biggl(\sup_{0\le t\le
T}\exp\biggl\{\int_t^Th(u)\,dW_u\biggr\}\biggr)^r\nonumber\\[-2pt]
&&\qquad=\EE\biggl(\exp\biggl\{\int_0^Th(u)\,dW_u\biggr\}\sup_{0\le t\le
T}\exp\biggl\{-\int_0^th(u)\,dW_u\biggr\}\biggr)^r\nonumber\\[-2pt]
&&\qquad\le\biggl(\EE\exp\biggl\{2r\int_0^Th(u)\,dW_u\biggr\}
\biggr)^{1/2}\\[-2pt]
&&\qquad\quad{}\times\biggl(\EE\sup_{0\le
t\le T}\exp\biggl\{-2r\int_0^th(u)\,dW_u\biggr\}\biggr)^
{1/2}\nonumber\\[-2pt]
&&\qquad=\exp\biggl\{r^2\int_0^Th(u)^2\,dr\biggr\}\Bigl(\EE\sup_{0\le
t\le
T}\eta_t^{2r}\Bigr)^{1/2}<\infty.\nonumber
\end{eqnarray}
For any $p^\prime\in(p,\frac{q}{2})$, by Doob's maximal inequality
and the H\"{o}lder inequality, (G3) and (\ref{ee.6.1}), we have
\begin{eqnarray*}
&&{\mathbb{E} \sup_{0\le i\le n}} |\delta Z_{t_i}^\pi|^p\\[-2pt]
&&\qquad\le C\mathbb{E} \biggl( \Bigl(\sup_{0\le\theta\le
T}|D_\theta\xi|\Bigr)^p\biggl(\sup_{0\le t\le
T}\exp\biggl\{\int_t^Th(r)\,dW_r\biggr\}\biggr)^p\\[-2pt]
&&\qquad\quad\hphantom{C\mathbb{E} \biggl(}
{}\times\biggl[\sup_{0\le
k\le n-1} \biggl|\int_{t_k}^{t_{k+1}} h(r) \,dW_r\biggr| \\[-2pt]
&&\qquad\quad\hspace*{41pt}{} + \sup_{0\le k\le n-1}\int_{t_k}^{t_{k+1}} |g(r)|\,dr
+\frac12
\sup_{0\le k\le n-1} \int_{t_k}^{t_{k+1}} h(r)^2 \,dr\biggr]^p %
\biggr)\\[-2pt]
&&\qquad\le C\biggl[\mathbb{E} \biggl( \Bigl(\sup_{0\le\theta\le
T}|D_\theta\xi|\Bigr)^{{pp^\prime}/({p^\prime-p})}\\[-2pt]
&&\hspace*{58.6pt}{}\times\biggl(\sup
_{0\le
t\le
T}\exp\biggl\{\int_t^Th(r)\,dW_r\biggr\}\biggr)^{{pp^\prime
}/({p^\prime-p})}\biggr)\biggr]^{({p^\prime-p})/{p^\prime}}\\[-2pt]
&&\qquad\quad{}\times\biggl[\EE\biggl(\sup_{0\le k\le n-1}
\biggl|\int_{t_k}^{t_{k+1}} h(r) \,dW_r\biggr| + \sup_{0\le k\le
n-1}\int_{t_k}^{t_{k+1}} |g(r)|\,dr \\[-2pt]
&&\qquad\quad\hspace*{136pt}{} +\frac12 \sup_{0\le k\le n-1}
\int_{t_k}^{t_{k+1}} h(r)^2 \,dr\biggr)^{p^\prime}
\biggr]^{{p}/{p^\prime}}\\[-2pt]
&&\qquad\le C\Bigl[\mathbb{E} \Bigl(\sup_{0\le\theta\le
T}|D_\theta\xi|\Bigr)^{{2pp^\prime}/({p^\prime-p})}
\Bigr]^{{p^\prime}/({2(p^\prime-p)})}\\[-2pt]
&&\qquad\quad{}\times\biggl[\EE\biggl(\sup_{0\le
t\le
T}\exp\biggl\{\int_t^Th(r)\,dW_r\biggr\}\biggr)^{{2pp^\prime
}/({p^\prime-p})}\biggr]^{{p^\prime}/({2(p^\prime-p)})}\\[-2pt]
&&\qquad\quad{}\times\biggl[\EE\sup_{0\le k\le n-1} \biggl|\int_{t_k}^{t_{k+1}}
h(r) \,dW_r\biggr|^{p^\prime} +\EE\sup_{0\le k\le
n-1}\biggl(\int_{t_k}^{t_{k+1}} |g(r)|\,dr\biggr)^{p^\prime}\\[-2pt]
&&\hspace*{148.3pt}\qquad\quad{} +\EE\sup_{0\le k\le n-1} \biggl(\int_{t_k}^{t_{k+1}} h(r)^2
\,dr\biggr)^{p^\prime} \biggr]^{{p}/{p^\prime}}\\[-2pt]
&&\qquad=C[I_1+I_2+I_3]^{{p}/{p^\prime}} .
\end{eqnarray*}

For any $r> 1$, by the H\"{o}lder inequality we can obtain
\begin{eqnarray*}
I_1 &=&\mathbb{E} \sup_{0\le k\le n-1} \biggl|\int_{t_k}^{t_{k+1}}
h(r)
\,dW_r\biggr|^{p^\prime}
\le\biggl\{ \mathbb{E} \sup_{0\le k\le
n-1}\biggl|\int_{t_k}^{t_{k+1}}
h(r) \,dW_r\biggr|^{p^\prime r}\biggr\}^{{1/r}} \\
&\le& \Biggl\{ \mathbb{E} \sum_{k=0}^{n-1}
\biggl|\int_{t_k}^{t_{k+1}} h(r) \,dW_r\biggr|^{p^\prime
r}\Biggr\}^{{1/r}} .
\end{eqnarray*}
For any centered Gaussian variable $X$, and any $\gamma\geq1$, we
know that %
\[
\mathbb{E} |X|^\gamma\le\tilde{C}^\gamma
\gamma^{{\gamma}/{2}} (\mathbb{E} |X|^2)^{{\gamma}/{2}},
\]
where $\tilde{C}$ is a constant independent of $\gamma$.
Thus, we can see that
\[
I_1 \le\Biggl( \tilde{C}^{p^\prime r} (p^\prime
r)^{{p^\prime r}/{2}} \sum_{i=0}^{n-1} \biggl(
\int_{t_i}^{t_{i+1}} h(r) ^2 \,dr \biggr)^{{p^\prime r}/{2}}
\Biggr)^{{1}/{r}}
\le C r^{{p^\prime}/{2}}
|\pi|^{{p^\prime}/{2}-{1/r}} .
\]
Take $r=\frac{2\log({1/|\pi|})}{p^\prime}$. Assume $%
|\pi|$ is small enough; then we have
\[
I_1\le C |\pi|^{{p^\prime}/{2}-{p^\prime}/({2\log
({1/|\pi|})})} \biggl( \log\frac{1%
}{|\pi|}\biggr)^{{p^\prime}/{2}} .
\]
It is easy to see that
\[
I_2 =\EE\sup_{0\le k\le n-1}\biggl(\int_{t_k}^{t_{k+1}}
|g(r)|\,dr\biggr)^{p^\prime} \le C |\pi|^{p^\prime}
\]
and%
\[
I_3=\EE\sup_{0\le k\le n-1} \biggl(\int_{t_k}^{t_{k+1}} h(r)^2
\,dr\biggr)^{p^\prime}\le C |\pi|^{p^\prime} .\vadjust{\goodbreak}
\]
Consequently, we obtain
%
%e5.8 ###
%
\begin{equation}\label{e.6.12}
\mathbb{E} \sup_{0\le i\le n} |\delta Z_{t_i}^\pi|^p \le C
|\pi|^{{p}/{2}-{p}/({2\log({1}/{|\pi|})})} \biggl( \log
\frac{1}{|\pi|}%
\biggr)^{{p/2}} .
\end{equation}
Applying recursively the scheme given by (\ref{e.6.5}), we
obtain
\[
Y_{t_i}^\pi= \mathbb{E} \Biggl( \xi+\sum_{k=i+1}^n f(t_k,
Y_{t_k}^\pi, Z_{t_k}^\pi) {\Delta}_{k-1}
\Big|\mathcal{F}_{t_i}\Biggr),\qquad i=n-1,n-2,\ldots,0.
\]
Therefore, for $i=n-1,n-2,\ldots,0$,
\[
|\delta Y_{t_i}^\pi|\le\mathbb{E} \Biggl( \sum_{k=i+1}^n
|f(t_k, Y_{t_k} , Z_{t_k} )-f(t_k, Y_{t_k}^\pi, Z_{t_k}^\pi)
|{\Delta}_{k-1} +|R_{t_i}^\pi|+|\delta\xi^\pi|
\Big|\mathcal{F}_{t_i}\Biggr) ,
\]
where $R_t^\pi$ is exactly the same as in Section \ref{sec3} and
$\delta
\xi^\pi=\xi-\xi=0$. In fact, we keep the term $\delta\xi^\pi$ to
indicate the role it plays as the terminal value.

For $j= n-1,n-2,\ldots,i$, we have
\begin{eqnarray*}
&&|\delta Y_{t_j}^\pi|\le\mathbb{E} \Biggl( \sum_{k=i+1}^n
|f(t_k, Y_{t_k} , Z_{t_k} )-f(t_k, Y_{t_k}^\pi, Z_{t_k}^\pi)
|{\Delta}_{k-1} \\[-3pt]
&&\hspace*{130.5pt}{}+\sup_{0\le t\le T}
|R_t^\pi|+|\delta\xi^\pi|\Big|\mathcal{F}_{t_j}\Biggr) .
\end{eqnarray*}
By Doob's maximal inequality and (\ref{e.6.12}), we obtain
\begin{eqnarray*} %\label{e.6.13}
&&\mathbb{E} \sup_{i\le j\le n} |\delta Y_{t_j}^\pi|^p \\[-3pt]
&&\qquad\le C\mathbb{E} \Biggl( \sum_{k=i+1}^n |f(t_k, Y_{t_k} ,
Z_{t_k} )-f(t_k, Y_{t_k}^\pi,
Z_{t_k}^\pi) |{\Delta}_{k-1} \Biggr)^p\\[-2pt]
&&\qquad\quad{} +C(|\pi|^{
{p/2}} +\EE|\delta\xi^\pi|^p)\\[-2pt]
&&\qquad\le C \Biggl\{ \mathbb{E} \Biggl( \sum_{k=i+1}^n | Y_{t_k} -
Y_{t_k}^\pi|{\Delta}_{k-1} \Biggr)^p +\mathbb{E} \Biggl(
\sum_{k=i+1}^n | Z_{t_k} - Z_{t_k}^\pi| {\Delta}_{k-1}
\Biggr)^p\Biggr\}\\[-2pt]
&&\qquad\quad{}+C(|\pi|^{{p/2}} +\EE|\delta\xi^\pi|^p)
\\[-2pt]
&&\qquad\le C_2 (T-t_i)^p {\mathbb{E} \sup_{i+1\le k\le n}} |
Y_{t_k}
- Y_{t_k}^\pi| ^p \\[-3pt]
&&\qquad\quad{} + C_3\biggl(|\pi|^{{p/2}-
{p}/({2\log({1}/{|\pi|}}))} \biggl( \log\frac{1}{|\pi|}%
\biggr)^{{p/2}}+\EE|\delta\xi^\pi|^p\biggr) ,
\end{eqnarray*}
where $C_2$ and $C_3$ are constants independent of the
partition~$\pi$.

We can obtain the estimate for ${\mathbb{E} \max_{0\le i\le n}}
|Y_{t_i}-Y_{t_i}^\pi|^p$ by using similar arguments to
analyze (\ref{e.5.10}) in Theorem \ref{t.5.2} to get the estimate
for ${\EE\sup_{0\le t\le T}}|Y_t-Y_t^\pi|$.\vadjust{\goodbreak}
\end{pf}

\section*{Acknowledgment}
We appreciate the referee's very constructive and
detailed comments to improve the presentation of this paper.

%suskaldyti doi

% imsref loaded by lrinkeviciute, 2011-03-18 15:02:15
%
% imsref loaded by lrinkeviciute, 2011-03-21 12:58:01

%
\printaddresses

\end{document}